\documentclass[aos,MSNbibl,seceqn,nameyear,dvips]{arximspdf}
\usepackage{accents}
\usepackage{graphicx}

\doi{10.1214/13-AOS1127} %kopijuoti is PTS
\volume{41}
\issue{1}
\pubyear{2013}
\firstpage{1780}
\lastpage{1815}

\makeatletter

\renewcommand{\underline}{\underaccent{\bar}}

\newtheorem{theorem}{Theorem}[section]
\newtheorem{prop}{Proposition}[section]
\newtheorem{lem}{Lemma}[section]

\newproclaim{defin}{Definition}[section]

\newcommand{\cA}{\mathcal{A}}
\newcommand{\cB}{\mathcal{B}}
\newcommand{\cE}{\mathcal{E}}
\newcommand{\cG}{\mathcal{G}}
\newcommand{\cN}{\mathcal{N}}
\newcommand{\cS}{\mathcal{S}}

\newcommand{\bP}{\mathbf{P}}

\newcommand{\bS}{\mathbf{S}}
\newcommand{\bone}{\mathbf{1}}
\newcommand{\Sh}{\hat{\Sigma}}
\newcommand{\N}{\mathcal{N}}
\newcommand{\Sy}{\mathbf{S}}
\newcommand{\R}{\mathbf{R}}

\newcommand{\Pro}{\mathbf{P}}

\newcommand{\E}{\mathbb{E}}

\newcommand{\DD}{\mathsf{D}}
\newcommand{\SDP}{\operatorname{\mathsf{SDP}}}
\newcommand{\MDP}{\operatorname{\mathsf{MDP}}}
\newcommand{\starDP}{*\mathsf{DP}}
\newcommand{\tr}{\operatorname{\mathbf{Tr}}}
\newcommand{\rk}{\operatorname{\mathbf{rank}}}

\newcommand{\ud}{\mathrm{d}}
\newcommand{\eps}{\varepsilon}

\makeatother

\begin{document}
\begin{frontmatter}

\title{Optimal detection of sparse principal components in high dimension}
\runtitle{Optimal detection of sparse principal components}

\begin{aug}
\author[A]{\fnms{Quentin} \snm{Berthet}\thanksref{t1}\ead[label=e1]{qberthet@princeton.edu}}
\and
\author[A]{\fnms{Philippe} \snm{Rigollet}\corref{}\thanksref{t2}\ead[label=e2]{rigollet@princeton.edu}}
\runauthor{Q. Berthet and P. Rigollet}
\affiliation{Princeton University}
\address[A]{Department of Operations Research\\
\quad and Financial Engineering\\
Princeton University\\
Princeton, New Jersey 08544\\
USA\\
\printead{e1}\\
\hphantom{E-mail: }\printead*{e2}} %adresu isvedimo komanda gale!
\end{aug}

\thankstext{t1}{Supported in part by a Gordon S. Wu fellowship.}
\thankstext{t2}{Supported in part by NSF Grants DMS-09-06424 and
CAREER-DMS-1053987.}

% HISTORY:
\received{\smonth{12} \syear{2012}}
\revised{\smonth{4} \syear{2013}}

% ABSTRACT
%
\begin{abstract}
We perform a finite sample analysis of the detection levels for sparse
principal components of a high-dimensional covariance matrix. Our
minimax optimal test is based on a sparse eigenvalue statistic. Alas,
computing this test is known to be NP-complete in general, and we
describe a computationally efficient alternative test using convex
relaxations. Our relaxation is also proved to detect sparse principal
components at near optimal detection levels, and it performs well on
simulated datasets. Moreover, using polynomial time reductions from
theoretical computer science, we bring significant evidence that our
results cannot be improved, thus revealing an inherent trade off
between statistical and computational performance.
\end{abstract}

% KEYWORDS
\begin{keyword}[class=AMS]
\kwd[Primary ]{62H25}
\kwd[; secondary ]{62F04}
\kwd{90C22}
\end{keyword}
\begin{keyword}
\kwd{High-dimensional detection}
\kwd{sparse principal component analysis}
\kwd{spiked covariance model}
\kwd{semidefinite relaxation}
\kwd{minimax lower bounds}
\kwd{planted clique}
\end{keyword}

\end{frontmatter}

%s1 #&#
\section{Introduction}
\label{SECintro}

The sparsity assumption has become preponderant in modern,
high-dimensional statistics. In the high dimension, low sample size
setting, where consistency seems to be hopeless, sparsity turns out to
be the statistician's salvation. It formalizes the a priori belief that
only a few parameters, among a large number of them, are significant
for the statistical task at hand. This paper explores a specific
high-dimensional problem, namely Principal Component Analysis (PCA).
Indeed, without further assumptions, classical PCA is known to produce
inconsistent estimators of the directions that explain the most
variance [\citet{JohLu09,Pau07,Nad08}]. For PCA, the \emph{spiked
covariance model} introduced by \citet{Joh01} provides a natural
setting for statistical problems. Namely, this model relies on the
assumption that there exists a small number of directions that explain
most of the variance. In this work, we assume that observations are
drawn from a multivariate Gaussian distribution with mean zero and
covariance matrix given by $I + \theta v v^\top$, where $I$ is the
identity matrix, $v$ is a unit norm sparse vector and $\theta>0$. Akin
to other models, the sparsity assumption drives both methods and
analysis in a wide variety of applications ranging from signal
processing to biology; see \citet
{AloBarNot99,Che11,JenOboBac10,WriGanYan11} for a few examples. Most
contributions to this problem
have focused on consistent estimation of the sparse principal component
$v$ under various performance measures; see, for example,
\citet{AmiWai09,Ma13,SheSheMar13,CaiMa12,VuLei12,BirJohNadPau12} and
the above references.

What if there is no sparse component? In other words, what if $\theta
=0$? From a detection standpoint, one may ask the following question:
How much variance should a sparse principal component explain in order
to be detectable by a statistical procedure? Answering this question
consists of (i) constructing a test that can detect this sparse
principal component when the associated variance is above a certain
level and (ii) proving that no test can detect such a principal
component below a certain level.

Optimal detection levels in a high-dimensional setup have recently
received a lot of attention. \citet
{DonJin04,IngTsyVer10,AriCanPla11,AriCanDur11} have studied the
detection of a sparse vector corrupted by noise under various sparsity
assumptions. More recently, this problem has been extended from vectors
to matrices by \citet{ButIng13}, Sun and Nobel
(\citeyear{SunNob08,SunNob13}) who propose to
detect a shifted sub-matrix planted in a Gaussian or binary random
matrix. While the notion of sub-matrix encodes a certain sparsity
structure, these two papers focus on the elementwise properties of
random matrices, unlike the blooming random matrix theory that focuses
on spectral aspects. \citet{AriBubLug12} studied a problem related to
sparse PCA detection, but closer to the shifted sub-matrix problem.
Their goal is to detect a shifted off-diagonal sub-matrix planted in a
covariance matrix. Their methods are not spectral either.

We extend the current work on detection in two directions. First, we
analyze detection in the framework of sparse PCA, and more precisely,
in the spiked covariance model. Second, we derive a finite sample
analysis of minimax optimality in this problem with results that hold
with high-probability, unlike most of the literature on detection where
an asymptotic framework is usually preferred. A notable exception is
the paper of \citet{AddBroDev10} where results are of the same flavor as
ours. Unlike the asymptotic analysis pioneered by \citet{Ing82} and
recently extended to sparse linear regression in \citet
{DonJin04,IngTsyVer10}, this finite sample analysis is not refined
enough to
exhibit a qualitative difference between testing and estimation.
Nevertheless, such results shed light on the delicate interplay between
the important parameters of the problem: ambient dimension, sample
size and sparsity.

The minimax optimal test statistic for our testing problem relies on
the so-called $k$-sparse largest eigenvalue of the empirical covariance
matrix. It captures the largest amount of empirical variance explained
by any $k$ of the original variables. It turns out that although this
statistic can be used to construct an optimal test, it raises
computational difficulties and can even be proved to be NP-complete in
general. As a result, a large body of the optimization literature on
this topic consists of numerical methods to overcome this issue; see,
for example, \citet{dAsBacEl-08,dAsGhaJor07,MaS11,LuZha11,JouNesRic10}
and references therein. Nevertheless, while these
numerical methods do produce a solution, their statistical properties
are rarely addressed for the estimation problem and never for the
detection problem. One of the approaches introduced by \citet
{dAsGhaJor07} uses a convexification technique called \emph
{semidefinite programming} (SDP). A major drawback of this technique is
that it may not output a sparse direction $\hat v$. Indeed,
semidefinite programs output matrices that are not rank-one in general,
and an ad hoc post-processing step is often required to turn this
matrix back into a unit vector. However, in the context of detection,
our goal is not to estimate the eigenvector $v$ but rather its
associated eigenvalue. This notable difference allows us to bypass SDP
optimization altogether, which is known to scale poorly in high
dimension. Inspired by the dual SDP formulation, we propose a simple
test procedure based on the \emph{minimum dual perturbation} (MDP)
that is easy to compute and for which we can derive near optimal
performance bounds for the detection problem. More importantly, we
bring supporting evidence to the tightness of the performance bounds
that we prove. Interestingly, this evidence builds upon a conjecture
from theoretical computer science. Indeed, a reduction to the \emph
{planted clique} problem shows that a better performance would
contradict a widely believed conjecture on the average-case complexity
of this problem. %We also provide a methodology to quickly find
%quantiles for our test statistic using Monte Carlo simulations.

Most of our analysis is performed in the model of sparse rank one
perturbation for the covariance matrix of Gaussian random vectors.
Nevertheless, our results are robust to variations around this model,
and we devote Section~\ref{SECweak} to discussing various weaker
assumptions under which our results still hold. In particular, our
results are more generally valid for sub-Gaussian observations and
weaker notions of sparsity. We also study the case where the distance
between the estimated and true covariance matrices is only controlled
in sup-norm, with high probability. This setup encompasses biased
estimators or adversarial noise.

The rest of the paper is organized as follows. In Section~\ref{SECpb},
we introduce the detection problem for sparse PCA. In Section
\ref{SECRMT}, we discuss various links with probabilistic results on
random matrix theory and more precisely, the asymptotic effect of a
principal component on the spectrum of a Wishart matrix. Minimax
detection levels are derived in Section~\ref{SECmain}, where in
particular, we introduce a test based on spectral methods and derive
the level at which it achieves detection of sparse principal components
with high probability. This level is proved to be optimal in a minimax
sense in Section~\ref{SEClow}. Unfortunately, this test cannot be
computed efficiently, and several relaxations are proposed in
Section~\ref{SECSDP}. For these convex methods, we derive suboptimal
levels that also hold under various weaker assumptions for which they
sometimes become optimal (Section~\ref{SECweak}). Moreover, using
arguments from computational complexity, we argue in
Section~\ref{SECCTLB} that even under the strongest assumptions of
this paper, these suboptimal levels are likely to be the best
achievable by the efficient relaxations. Specifically, we show that
proving better bounds for these methods would lead to a contradiction
of the hidden clique conjecture, which is widely believed to be true.
The numerical performance of our test and in particular its
suboptimality, is illustrated in Section~\ref{SECsim}.\vspace*{9pt}

\textsc{Notation.} The space of $d \times d$ symmetric real matrices
is denoted by $\Sy_d$. We write $Z \succeq0$ whenever $Z$ is
semidefinite positive.

The elements of a vector $v \in\R^d$ are denoted by $v_1,\ldots, v_d$
and similarly, a matrix $Z$ has element $Z_{ij}$ on its $i$th row and
$j$th column. For any $q > 0$, $|v|_q$ denotes the $\ell_q$ ``norm'' of
a vector $v$ and is defined by $|v|_q=(\sum_{j}|v_j|^q)^{1/q}$.
Moreover, we denote by $|v|_0$ its so-called $\ell_0$ ``norm,'' that
is, its number of nonzero elements. Furthermore, by extension, for $Z
\in\Sy_d$, we denote by $|Z|_q$ the $\ell_q$ norm of the vector formed
by the entries of $Z$. We also define for $q\in[0,2)$ the set
$\cB_q(R)$ of unit vectors within the $\ell_q$-ball of radius $R>0$
\[
\cB_q(R) = \bigl\{v \in\R^p\dvtx |v|_2=1,
|v|_q \le R\bigr\}.
\]

The trace and rank functionals are denoted by $\tr$ and $\rk$,
respectively, and have their usual definition. The identity matrix in
$\R^d$ is denoted by $I_d$. For a finite set $S$, we denote by $|S|$
its cardinality. We also write $A_S$ for the $|S|\times|S|$ submatrix
with elements $(A_{ij})_{i,j \in S}$, and $v_S$ for the vector of $\R
^{|S|}$ with elements $v_i$ for \mbox{$i \in S$}. Finally, for two real
numbers $a$ and $b$, we write $a \wedge b=\min(a,b)$, $a \vee b = \max
(a,b)$ and $a_+ = a \vee0 $.

%s2 #&#
\section{Statement of the hypothesis testing problem}
\label{SECpb}

Let $X_1,\ldots,X_n$ be $n$ i.i.d. copies of a random variable $X$ in
$\R^p$. Our objective is to perform the following test:
\begin{eqnarray*}
H_0 \dvtx X &\sim& \N(0,I_p),
\\
H_1 \dvtx X &\sim& \N\bigl(0,I_p+\theta vv^\top
\bigr),\qquad v \in\cB_0(k),
\end{eqnarray*}
where $\theta>0$, and we remind the reader that $\cB_0(k)$ is the set
of $k$-sparse unit vectors.
Note that the model under $H_1$ is an adaptation of the spiked
covariance model since\vadjust{\goodbreak} it only allows $v$ to be $k$-sparse on the unit
Euclidean sphere. This is precisely the model of sparse PCA introduced
in \citet{JohLu09}. In particular, the distribution of $X$ under $H_1$
is invariant under rotation of the $k$ relevant variables. We use this
simplified model for reasons of clarity: to highlight the importance of
relative variance, only one direction $v$ is used for signal, and only
one parameter $\theta$ is used to express the signal-to-noise ratio.
Note that our upper and lower bounds for optimal testing are valid for
the general hypotheses
\begin{eqnarray*}
H_0 \dvtx X &\sim& \N(0,\Sigma_0),\qquad \lambda_{\max}^k(
\Sigma_0) \leq1,
\\
H_1 \dvtx X &\sim& \N(0,\Sigma_1),\qquad \lambda_{\max}^k(
\Sigma_1) \geq1+\theta,
\end{eqnarray*}
where $\lambda_{\max}^k$, is the $k$-sparse eigenvalue defined
in (\ref
{EQlambdakmax}) below. In particular, the model under $H_1$ encompasses
that of \citet{AmiWai09}.

Let $\Sigma=\E[XX^\top]$ denote the covariance matrix of the
centered random vector~$X$, and denote by $\hat{\Sigma}$ the
empirical covariance matrix defined by
%
%e2.1 #&#
%
\begin{equation}
\label{EQdefEmpCov} \Sh= \frac{1}{n} \sum_{i=1}^n
X_i X_i^\top.
\end{equation}

We say that a test \emph{discriminates} between $H_0$ and $H_1$ with
probability $1-\delta$ if both type I and type II errors have a
probability smaller than $\delta$. Our goal is therefore to find a
statistic $\varphi(\hat{\Sigma})$ and quantiles $\tau_0 < \tau_1$,
depending on ($p,n,k,\delta$) such that
\[
\Pro_{H_0}\bigl(\varphi(\hat{\Sigma})> \tau_0\bigr) \leq
\delta,\qquad \Pro_{H_1}\bigl(\varphi(\hat{\Sigma})< \tau_1\bigr)
\leq\delta.
\]
For $\tau\in[\tau_0, \tau_1]$ define the test
\[
\psi(\hat{\Sigma}) = \bone{\bigl\{ \varphi(\hat{\Sigma})> \tau\bigr\}},
\]
where $\bone\{\cdot\}$ denotes the indicator function.
As desired, this test discriminates between the two hypotheses with
probability $1-\delta$. We assume that the user is testing for a
specific sparsity $k$. Nevertheless, using a Bonferroni correction,
this test can be performed for various values of $k$ if needed.

Note that throughout the paper, we assume that all of the parameters
$(k,n,p)$ are known so that $\tau_0$ and $\tau_1$ are easily determined.

%s3 #&#
\section{Link with random matrix theory}
\label{SECRMT}

%Note that under the null hypothesis, the sample covariance matrix $
%such matrices has been extensively studied and is fairly well
%understood. We give below a quick overview of the results that are
%relevant to our problem.

%s3.1 #&#
\subsection{Spectral methods}

It is not hard to see that, under $H_1$, for any $\theta>0$, $v$~is an
eigenvector associated to the largest eigenvalue of the population
covariance matrix $\Sigma$. Moreover, if $\hat\Sigma$ is close to
$\Sigma$ in spectral norm, then its largest eigenvector should be a
good candidate to approximate $v$. It is therefore natural to consider
spectral methods for the spiked covariance model. Understanding the
behavior of our test statistic under both the null and the alternative
is key in proving that it discriminates between the hypotheses.

Spectral convergence of the empirical covariance matrix to the true
covariance matrix has received some attention recently [see, e.g.,
\citet{CaiZhaZho10,BicLev08,El-08}] under various elementwise
sparsity assumptions and using thresholding methods. However, since our
assumption allows for relevant variables to produce arbitrarily small
entries under the alternative hypothesis, we cannot use such results.
A~natural statistic to discriminate between the null and the alternative
would be, for example, the largest eigenvalue of the covariance matrix.

Spectral properties of random matrices have received a lot of attention
from both a statistical and probabilistic perspective. We devote the
rest of this section to reviewing some of the classical results from
random matrix theory, and we argue that even in moderate dimension, the
largest eigenvalue cannot discriminate between the null and alternative
hypotheses.

It is easily seen that for any unit vector $v$,
%
%e3.1 #&#
%
\begin{equation}
\label{EQlamax} \lambda_{\max}(I_p) = 1 \quad\mbox{and}\quad
\lambda_{\max}\bigl(I_p + \theta vv^\top\bigr) = 1 +
\theta.
\end{equation}

If we could allow, for a fixed $p$, to let $n$ go to infinity, the
consistency of the estimator $\hat{\Sigma}$ (for fixed $p$, entry by
entry) and the continuity of the largest eigenvalue as a function of
the entries of a matrix would imply that the largest eigenvalue can be
used to discriminate between the two alternatives, at least asymptotically.

However, in a high-dimension setting, where $p$ is typically much
larger than $n$, the behavior of $\lambda_{\max}(\hat{\Sigma})$
under the null hypothesis is quite different.
If $p/n \rightarrow\alpha>0$, \citet{Gem80} showed that, in
accordance with the Marcenko--Pastur distribution, we have
\[
\lambda_{\max}(\hat{\Sigma}) \rightarrow(1+\sqrt{\alpha})^2
> 1,
\]
where the convergence holds almost surely see also \citet{Joh01,Bai99}
and references therein. Moreover, \citet{YinBaiKri88}
established that finite fourth moment is a necessary and sufficient
condition for this almost sure convergence to hold.
Furthermore, since $\hat{\Sigma} \succeq0$, its number of positive
eigenvalues is equal to its rank (which is at most $n$), and we have
\[
\lambda_{\max}(\hat{\Sigma}) \geq\frac{1}{\rk(\hat{\Sigma
})}\sum
_{i=1}^p{\lambda_i(\hat{\Sigma}}) \geq
\frac{1}{n} \tr(\hat{\Sigma}) = \frac pn \frac{\sum_{i=1
}^n|X_i|_2^2}{np}.
\]
Note that under $H_0$, it holds $\sum_{i=1
}^n|X_i|_2^2 \sim\chi^2_{np}$. Hence almost surely, for $p/n
\rightarrow\infty$, we have $\lambda_{\max}(\hat{\Sigma})
\rightarrow\infty$.

These two results hint at an intrinsic limitation of the largest
eigenvalue statistic: its fluctuations are too large to discriminate
between the two hypotheses in a ``large $p$/small $n$'' scenario unless
the signal strength $\theta$ is very strong.

In the next subsection, we show that the above argument can be made
formal using spectral results in random matrix theory.

%s3.2 #&#
\subsection{Finite rank perturbations of covariance matrices}
In a moderate-dimen\-sional regime, where $p/n \rightarrow\alpha\in
(0,1)$, \citet{BaiAroPec05} describe a phase transition for the
spectral behavior of the sample covariance matrix $\Sh$ of complex
Gaussian vectors between two different regimes. This phenomenon is now
widely known as the BBP transition for the name of the authors. The
same phenomenon, for real random variables, was subsequently
established in \citet{BaiSil06}.

Qualitatively,\vspace*{1pt} there exists a critical value $\theta^*$ such that if
$\theta>\theta^*$, the spectrum of $\Sh$ exhibits an isolated
eigenvalue significantly larger than the others, and such that if
$\theta<\theta^*$, the spectrum has a similar behavior under the two
hypotheses. More precisely, Theorem 1.1 of \citet{BaiSil06} implies
that under $H_1$, the largest eigenvalue will either exhibit an
important concentration around a deterministic value strictly larger
than $1+\theta$ if the perturbation is strong enough, or around the
upper edge of the Marcenko--Pastur distribution, as if the perturbation
was nonexistent, when it is too weak. The critical level is $\theta^*
= \sqrt{\alpha}$, and suggests a minimum signal level of order $\sqrt
{p/n}$ which is high already when $p$ is of the order of $n$.

These results are even proved to hold for weakened assumptions on the
distribution of the vectors, in \citet{FerPec09}. On the statistical
side, these are coherent with the detection levels shown in \citet
{OnaMorHal11} for testing of the sphericity hypothesis with no
assumption on the alternative, by spectral methods.

%s4 #&#
\section{Sparse principal component detection}
\label{SECmain}

In sparse principal component detection, we are testing the existence
of a sparse direction $v$ with a significantly higher explained
variance $v^\top\Sigma v$ than any other direction. To exploit the
sparsity assumption, we use the fact that only a small submatrix of the
covariance is affected by the perturbation. Let $A$ be a $p \times p$
matrix and fix $k <p$. We define the $k$-sparse largest\setcounter{footnote}{2}\footnote{In
the rest of the paper, we drop the qualification ``largest'' since we
only refer to this one.} eigenvalue by
%
%e4.1 #&#
%
\begin{equation}
\label{EQlambdakmax} \lambda^k_{\max}(A) = \max
_{|S| = k} \lambda_{\max}(A_S).
\end{equation}
It can be defined equivalently to (\ref{EQlambdakmax}) by
%
%e4.2 #&#
%
\begin{equation}
\label{EQlambdakmax2} \lambda^k_{\max}(A) = \max
_{x \in\cB_0(k)} x^\top A x.
\end{equation}

Therefore, we study the behavior of the test statistic $\varphi(\hat
{\Sigma}) = \lambda^k_{\max}(\hat{\Sigma})$ under both hypotheses.

%s4.1 #&#
\subsection{Deviation bounds for the $k$-sparse eigenvalue}

Optimal detection levels are governed by the deviations of the test
statistic $ \lambda^k_{\max}(\hat{\Sigma})$ both under the null and
the alternative hypotheses.
We begin with the following proposition, which guarantees that our test
statistic remains large enough under the alternative hypothesis.

%pr4.1 #&#
%
\begin{prop}
\label{THH1lk}
Under $H_1$, we have with probability $1-\delta$,
\[
\lambda^k_{\max}(\hat{\Sigma}) \geq1+ \theta- 2(1+ \theta)
\sqrt{\frac{\log(1/\delta)}{n}}.
\]
\end{prop}
\begin{pf}
Under $H_1$, there exists a unit vector $v$ with sparsity $k$, such
that $X \sim\N(0,I_p+\theta vv^\top)$. Therefore, we have
\[
\lambda^k_{\max}(\hat{\Sigma}) \geq v^\top\hat{
\Sigma} v = \frac{1}{n} \sum_{i=1}^{n}
\bigl(X_i^\top v\bigr)^2
\]
by definition of $\hat{\Sigma}$. Since $X \sim\N(0,I_p+\theta
vv^\top)$, we have $X^\top v \sim\N(0,1+\theta)$.

Define the random variable
\[
Y = \frac{1}{n} \sum_{i=1}^{n}
\biggl(\frac{(X_i^\top v)^2}{1+\theta
} - 1 \biggr).
\]
Using Lemma~\ref{LEMLM}, we get for any $t>0$, that
\[
\mathbf{P} (Y \leq-2\sqrt{t/n} ) \leq e^{-t}.
\]
Hence, taking $t = \log(1/\delta)$ yields the desired inequality.
\end{pf}

Note that our proof relies only on the existence of a sparse vector $v$
associated to the eigenvalue $(1+\theta)$ of the population covariance
matrix $\Sigma$. In particular, the result of Proposition
\ref{THH1lk} extends to more general alternative hypotheses, as long as
they satisfy this condition.

Note that much more than detection can actually be achieved under this
model. Indeed, \citet{AmiWai09} prove optimal rates of support recovery
when $\theta$ is known and large enough, and for $v$ taking only
values in $\{0, \pm1/\sqrt k\}$.

We now study the behavior of the $k$-sparse eigenvalue under the null
hypothesis, that is, for a Wishart matrix with mean $I_p$. We adapt a
technique from \citet{Ver12} to obtain the desired deviation bounds.

%pr4.2 #&#
%
\begin{prop}
\label{THH0lk}
Under $H_0$, with probability $1-\delta$
\[
\lambda_{\max}^k(\hat{\Sigma}) \leq1+ 4 \sqrt{
\frac{k \log(9ep
/ k)+\log(1/\delta)}{n}} + 4 \frac{k \log(9ep / k)+\log(1/\delta
)}{n}.
\]
\end{prop}
%
%We need a few lemmas to prove this result. The following lemma can be
%found in \citet{Ver12}. It relies on the construction of a $1/4$ net
%over the unit sphere of $\R^k$.
% A well-known approach is the $\varepsilon$-net on the unit sphere,
%here with $\varepsilon= 1/4$, found in \citet{VershTut}.
%
\begin{pf}
Using a $1/4$-net over the unit sphere of $\R^k$, it can be easily
shown [see, e.g., \citet{Ver12}] that there exists a subset $\mathcal
{N}_k$ of the unit sphere of $\R^k$, with cardinality smaller than
$9^k$, such that for any $A \succeq0$
%
%e4.3 #&#
%
\begin{equation}
\label{EQepsnet} \lambda_{\max}(A) \leq2 \max_{x \in\mathcal{N}_k}
x^\top A x.
\end{equation}
Under $H_0$, since $\hat\Sigma$ is positive semidefinite, we have
\[
\lambda_{\max}^k(\hat{\Sigma}) = 1+ \max
_{|S| = k} \bigl\{\lambda_{\max}(\hat{\Sigma}_S)-1
\bigr\}.
\]
%
%where the maximum in the right-hand side is taken over all subsets of $
For all $u\in\R^k, |u|_2=1$ and $S\subset\{1,\ldots, p\}$ such
that $|S|=k$, let $\tilde{u} \in\R^p$ be the vector with support in
$S$ such that $\tilde{u}_S = u$. We have
\[
u^\top\hat{\Sigma}_S u - 1= \tilde{u}^\top\Sh
\tilde{u} -1= \frac{1}{n} \sum_{i =1}^n
\bigl[\bigl(\tilde{u}^\top X_i\bigr)^2 -1
\bigr].
\]
Since $|\tilde{u}|_2=|u|_2=1$, Lemma~\ref{LEMLM} yields for any $t>0$,
%
%e4.4 #&#
%
\begin{equation}
\label{EQchi1} \mathbf{P} \Biggl(\frac{1}{n} \sum
_{i = 1}^n \bigl[\bigl(\tilde{u}^\top
X_i\bigr)^2 -1 \bigr] \geq2\sqrt{\frac{t}{n}} + 2
\frac{t}{n} \Biggr) \leq e^{-t}.
\end{equation}

For any $S \subset\{1,\ldots, p\}$, define $\R^S$ to be the subset
of $\R^p$ defined such that $x \in\R^S$ if and only if $x_j=0,
\forall j \notin S$. Let $\mathcal{N}_{k}(S)$ be a subset of the
unit sphere of~$\R^S$, with cardinality smaller than $9^k$ such that
for any $A \succeq0$, inequality (\ref{EQepsnet}) holds with $\cN
_k=\cN_k(S)$. Fix $t>0$ and define the event $\mathcal{A}_S$ by
\[
\mathcal{A}_S = \biggl\{ \lambda_{\max}(\hat{
\Sigma}_S)-1 \geq4\sqrt{\frac{t}{n}} + 4 \frac{t}{n}
\biggr\}.
\]
Observe that a union bound over the elements of $\cN_k(S)$ together
with (\ref{EQchi1}) yields that for any $t>0$,
\[
\mathbf{P} (\mathcal{A}_S ) \le\mathbf{P} \Biggl( \max
_{v \in
\mathcal{N}_{k}(S)}\frac{1}{n} \sum_{i=1}^n
\bigl(v^\top X_i \bigr)^2 -1\geq2\sqrt{
\frac{t}{n}} + 2 \frac{t}{n} \Biggr)\le9^k
e^{-t}.
\]
Let now $\cA$ be the event defined by
\[
\mathcal{A} = \bigcup_{|S|=k} \mathcal{A}_S
= \biggl\{ \max_{|S| =
k} \bigl\{\lambda_{\max}(\hat{
\Sigma}_S)-1 \bigr\} \geq4\sqrt{\frac{t}{n}} + 4
\frac{t}{n} \biggr\}.
\]
Therefore, by a union bound on the ${p \choose k}$ subsets $S$ of $\{
1,\ldots, p\}$ that have cardinality~$k$, we get
\[
\mathbf{P} \biggl( \lambda_{\max}^k(\hat{\Sigma}) \geq1+ 4
\sqrt{\frac{t}{n}} + 4 \frac{t}{n} \biggr) = \mathbf{P}(\mathcal{A})
\le\pmatrix{p \cr k} 9^k e^{-t}.
\]
To complete our proof, it is sufficient to use the standard inequality
${p \choose k} \leq(\frac{ep}{k} )^k$ and to take $t = k
\log(9ep / k)+\log(1/\delta)$.
\end{pf}

%s4.2 #&#
\subsection{\texorpdfstring{Hypothesis testing with $\lambda^k_{\max}$}
{Hypothesis testing with lambda k max}}

Using these results, we have, with the notation from Section~\ref{SECpb},
\[
\Pro_{H_0}\bigl(\lambda^k_{\max}(\hat{\Sigma})>
\tau_0\bigr) \leq\delta,\qquad \Pro_{H_1}\bigl(
\lambda^k_{\max}(\hat{\Sigma})< \tau_1\bigr) \leq
\delta,
\]
where $\tau_0$ and $\tau_1$ are given by
\begin{eqnarray*}
\tau_0 &=& 1+ 4 \sqrt{\frac{k \log(9ep / k)+\log(1/\delta)}{n}} + 4
\frac{k \log(9ep /
k)+\log(1/\delta)}{n},
\\
\tau_1 &=& 1+ \theta- 2(1+ \theta)\sqrt{\frac{\log({1/\delta})}{n}}.
\end{eqnarray*}

Whenever $\tau_1>\tau_0$, we take $\tau\in[\tau_0,\tau_1]$ and
define the following test:
\[
\psi(\hat{\Sigma}) = \bone{\bigl\{ \lambda^k_{\max}(\hat{
\Sigma})> \tau\bigr\}}.
\]
It follows from the previous subsection that the test discriminates
between $H_1$ and $H_0$ with probability $1-\delta$.
It remains to find for which values of $\theta$ the condition $\tau
_1>\tau_0$ holds. It corresponds to our minimum detection level.
%
%th4.1 #&#
%
\begin{theorem}
\label{CORbartheta}
Assume that $k, p, n$ and $\delta$ are such that $\bar\theta\leq1$, where
%
%e4.5 #&#
%
\begin{eqnarray}
\label{EQbartheta}
\bar\theta&:=&4 \sqrt{\frac{k \log({9ep}/{k})+\log
(1/{\delta})}{n}} + 4 \frac{k \log({9ep}/{k})+\log(1/{\delta})}{n} \nonumber\\[-8pt]\\[-8pt]
&&{}+
4\sqrt{\frac{\log(1/{\delta})}{n}}.\nonumber
\end{eqnarray}
Then, for any $\theta> \bar\theta$ and for any $\tau\in[\tau
_0,\tau_1]$, the test $\psi(\hat{\Sigma}) = \bone\{ \lambda
^k_{\max}(\hat{\Sigma})> \tau\}$ discriminates between $H_0$ and
$H_1$ with probability $1-\delta$.
\end{theorem}
If we consider high-dimensional asymptotic regimes, for large $p,n,k$,
taking $\delta= p^{-\beta}$ with $\beta>0$, provides a sequence of
tests $\psi_n$ that discriminate between $H_0$ and $H_1$ with
probability converging to 1, for any fixed $\theta>0$, as soon as $k
\log(p) /n \rightarrow0$.

%s5 #&#
\section{Minimax lower bounds for detection}
\label{SEClow}

The goal of this section is to prove that for any $\nu>0$, there
exists $\underline{\theta}_\nu$ (\ref{EQutheta}) of the same order
as $\bar\theta$ (up to logarithmic terms), and such that if $\theta
<\underline{\theta}_\nu$, then no test can discriminate between
$H_0$ and $H_1$ with probability greater than $\frac12 + \nu$.
%We will see that this result can be achieved up to logarithmic terms
%that vanish for interesting regimes of $p, n$ and $k$.
%Throughout this section, we assume that $\theta<1/\sqrt{2}$.
Recall that $\Pro^{ n}$ denotes the joint distribution of $n$ i.i.d.
random variables with distribution~$\Pro$.

%th5.1 #&#
%
\begin{theorem}
\label{THinfminimax}
Fix $\nu>0$. There exists a constant $C_\nu>0$ defined in
(\ref{EQdefCnu}) such that if
%
%e5.1 #&#
%
\begin{equation}
\label{EQutheta} \theta< \underline{\theta}_\nu:=\sqrt{
\frac{ k\log(C_\nu
p/k^2+ 1 )}{n}} \wedge\frac{1}{\sqrt{2}},
\end{equation}
it holds
\[
\inf_{\psi} \Bigl\{ \Pro_{0}^{ n}(\psi=1)
\vee\max_{v \in\cB
_0(k)} \Pro^{ n}_v(\psi=0) \Bigr
\} \geq\frac12 - \nu,
\]
where the infimum is taken over all possible tests, that is, measurable
functions of the $n$ observations, that take values in $\{0,1\}$.
\end{theorem}

In order to find lower bounds for the probability of error, we study
the $\chi^2$ distance between probability measures; see, for
example, \citet{Tsy09}, Chapter~2. For any $v \in\R^p$ such that
$|v|_2=1$, define the matrix $\Sigma_v=I_p+\theta vv^\top$, and let
$\Pro_v$ denote the distribution of a Gaussian random variable $X \sim
\cN(0, \Sigma_v)$. Moreover, let $\cS=\{S \subset\{1,\ldots,p\}\dvtx
|S| = k\}$, and for any $S \in\cS$, define $u(S) \in\R^p$ to be
the unit vector with $j$th coordinate equal to $1/\sqrt{k}$ if $j \in
S$ and $0$ otherwise. Finally, define the Gaussian mixture $\Pro
_{\mathcal{S}}$ by
\[
\Pro_{\mathcal{S}} = \frac{1}{|\mathcal{S}|} \sum_{S \in\mathcal
{S}}
\Pro_{u(S)}.
\]

We write for simplicity $\Pro_{S}=\Pro_{u(S)}$ when this leads to no
confusion. Our proof relies on the following lemma.
%
%le5.1 #&#
%
\begin{lem}
\label{LEMcross}
For any $S,T \in\mathcal{S}$ and any $\theta<1$, it holds
\[
\E_{\Pro_{0}} \biggl(\frac{\ud\Pro_{S}}{\ud\Pro_{0}}\,\frac{\ud
\Pro_{T}}{\ud\Pro_{0}} \biggr) = \bigl(1-
\theta^2 \bigl(u(S)^\top u(T)\bigr)^2
\bigr)^{-1/2}.
\]
\end{lem}
\begin{pf}
Fix $S \in\cS$, and observe that
\[
\frac{\ud\Pro_{S}}{\ud\Pro_{0}}(X) = \frac{\det(I_p)^{1/2}}{\det
(\Sigma_{u(S)})^{1/2}} \frac{\exp(-X^\top\Sigma
_{u(S)}^{-1} X/2)}{\exp(-X^\top I_p^{-1} X/2)}.
\]
Furthermore, since $\det(I_p) = 1$ and $|u(S)|_2=1$, we get by
Sylvester's determinant theorem that
\[
\det(\Sigma_{u(S)}) = \det\bigl(I_p + \theta u(S)
u(S)^\top\bigr) = \det\bigl(I_1 + \theta
u(S)^\top u(S)\bigr) =1+\theta.
\]
Moreover, the Sherman--Morrison formula yields
\[
\Sigma_{u(S)}^{-1} = \bigl(I_p + \theta u(S)
u(S)^\top\bigr)^{-1} = I_p- \frac
{\theta u(S) u(S)^\top}{1+\theta}.
\]
By substitution, the above three displays yield
\[
\frac{\ud\Pro_{S}}{\ud\Pro_{0}}(X) = \frac{1}{\sqrt{1+\theta}} \exp
\biggl(\frac{1}{2}
\frac{\theta}{1+\theta} \bigl(X^\top u(S)\bigr)^2 \biggr)
\]
and
%
%e5.2 #&#
%
\begin{equation}
\label{EQcross} \frac{\ud\Pro_{S}}{\ud\Pro_{0}}\frac{\ud\Pro_{T}}{\ud
\Pro
_{0}}(X) =\frac{1}{1+\theta}
\exp\bigl(X^\top M X\bigr),
\end{equation}
where $M$ is defined by
\[
M:= \frac{1}{2} \frac{\theta}{1+\theta} \bigl(u(S) u(S) ^\top+ u(T)
u(T)^\top\bigr).
\]
Note that $M$ has at most two nonzero eigenvalues given by
\[
\lambda_1 = \frac{1}{2} \frac{\theta}{1+\theta}
\bigl(1+u(S)^\top u(T)\bigr)<\frac1{2} \quad\mbox{and}\quad
\lambda_2 = \frac{1}{2} \frac{\theta}{1+\theta}\bigl(1-u(S)^\top
u(T)\bigr)<\frac12,
\]
and let $\Lambda$ denote the diagonal matrix with elements $(\lambda
_1,\lambda_2,0,\ldots,0) \in\R^p$.
Together with (\ref{EQcross}), it yields
\begin{eqnarray*}
\E_{\Pro_{0}} \biggl(\frac{\ud\Pro_{S}}{\ud\Pro_{0}}\frac{\ud
\Pro_{T}}{\ud\Pro_{0}} \biggr) &=&
\frac{1}{1+\theta} \E_{\Pro_{0}} \bigl[\exp\bigl(X^\top M X\bigr)
\bigr]
\\
&=& \frac{1}{1+\theta} \E_{\Pro_{0}} \bigl[\exp\bigl(X^\top\Lambda
X\bigr)\bigr]
\\
&=& \frac{1}{1+\theta} \E_{\Pro_{0}}\bigl[\exp\bigl(\lambda_1
X_1^2\bigr)\bigr] \E_{\Pro_{0}}\bigl[\exp\bigl(
\lambda_2 X_2^2\bigr)\bigr]
\\
&=& \frac{1}{1+\theta} \bigl[(1-2\lambda_1) (1-2
\lambda_2) \bigr]^{-1/2},
\end{eqnarray*}
where, in the second equality, the substitution of $M$ by $\Lambda$ is
valid by rotational invariance of the distribution of $X$ under $\Pro
_{0}$. The last equation yields the desired result.
\end{pf}

We now turn to the proof of Theorem~\ref{THinfminimax}.

\begin{pf*}{Proof of Theorem~\ref{THinfminimax}}
Observe now that
\[
\chi^2(\Pro_{\mathcal{S}},\Pro_{0}) = \E_{\Pro_{0}}
\biggl[ \biggl( \frac{\ud\Pro_{\mathcal{S}}}{\ud\Pro_{0}}-1 \biggr)^2
\biggr] =
\frac{1}{|\mathcal{S}|^2} \sum_{S,T \in\mathcal{S}} \E_{\Pro
_{0}}
\biggl(\frac{\ud\Pro_{S}}{\ud\Pro_{0}}\frac{\ud\Pro
_{T}}{\ud\Pro_{0}} \biggr) -1.
\]
Lemma~\ref{LEMcross} together with the fact $u(S)^\top u(T) = |S
\cap T|/k$ yield
\[
\chi^2(\Pro_{\mathcal{S}},\Pro_{0}) = \sum
_{r=0}^k \biggl\{ \frac
{\mathcal{C}(\mathcal{S},r)}{|\mathcal{S}|^2} \biggl(1-
\theta^2 \frac{r^2}{k^2} \biggr)^{-1/2} \biggr\}-1,
\]
where $\mathcal{C}(\mathcal{S},r)$ denotes the number of subsets $S,T
\in\cS$ such that $|S \cap T|=r$.
Let $S, T$ be chosen uniformly at random in $\mathcal{S}$, and observe
that $\bP(|S\cap T|=r)=\bP(R=r)$, where $R=|S \cap\{1,\ldots,k\}|$.
Jensen's inequality yields
\begin{eqnarray*}
\chi^2\bigl(\Pro_{\mathcal{S}}^{ n},
\Pro_{0}^{ n}\bigr) &=& \prod_{i=1}^n
\bigl(1+\chi^2(\Pro_{\mathcal{S}},\Pro_{0}) \bigr)-1
\\
&\le&\E_{S,T} \biggl\{ \biggl[ 1- \theta^2
\frac{|S\cap
T|^2}{k^2} \biggr]^{-n/2} \biggr\} -1\\
&=& \E_{R} \biggl\{
\biggl[ 1- \theta^2 \frac{R^2}{k^2} \biggr]^{-n/2} \biggr
\}-1,
\end{eqnarray*}
where $\E_{S,T}$ denotes the expectation with respect to the random
subsets $S,T$ and $\E_R$ the expectation with respect to $R$.

Using now the convexity inequality $(1-t)^{-n/2} \leq e^{
{nt}/({2(1-t)})} \leq e^{nt}$ valid for $1-t \geq1/2$, and noticing that
$R \le k$, the above display leads to
%
%e5.3 #&#
%
\begin{equation}
\label{EQchiesp} \chi^2\bigl(\Pro_{\mathcal{S}}^{ n},
\Pro_{0}^{ n}\bigr) \le\E_{R} \biggl[ \exp
\biggl( \frac{ n \theta^2R}{k} \biggr) \biggr]-1.
\end{equation}

Define $\mu^2= n\theta^2/k$. We have, as in \citet
{AddBroDev10,AriBubLug12}, that
\begin{eqnarray*}
\E_R \bigl[e^{\mu^2 R} \bigr] &=& \E_S \Biggl[\prod
_{i=1}^k \exp\bigl(\mu^2
\bone{\{i\in S\}}\bigr) \Biggr]
\\
&\leq&\prod_{i=1}^k\E_S \bigl[
\exp\bigl(\mu^2 \bone{\{i\in S\}}\bigr) \bigr]\leq\biggl(
\bigl(e^{\mu^2}-1 \bigr)\frac{k}{p}+1 \biggr)^k.
\end{eqnarray*}
The first inequality holds by the negative association [see, e.g.,
\citet{AddBroDev10}, Section 3] of negatively correlated dependent
random variables. Assume now that $\theta< \underline{\theta}_\nu$.
It yields
\[
\biggl( \bigl(e^{\mu^2}-1 \bigr)\frac{k}{p}+1
\biggr)^k\le\biggl( \biggl(\frac{C_\nu p}{k^2} \biggr)
\frac{k}{p} +1 \biggr)^k\le\biggl(1+\frac
{C_\nu}{k}
\biggr)^k \le e^{C_\nu}.
\]
Together with (\ref{EQchiesp}), the previous two displays yield
%
%e5.4 #&#
%
\begin{equation}
\label{EQchitot} \chi^2\bigl(\Pro_{\mathcal{S}}^{ n},
\Pro_{0}^{ n}\bigr) \le e^{C_\nu
}-1.
\end{equation}
We are now in a position to apply standard results from minimax theory. Define
%
%e5.5 #&#
%
\begin{equation}
\label{EQdefCnu} C_\nu:= \log\biggl[\bigl(1+8\nu^2
\bigr)\wedge\log\biggl(\frac{e}{2-4\nu} \biggr) \biggr],
\end{equation}
and note that for all measurable tests $\psi$, we have
\begin{eqnarray*}
\Pro_{0}^{ n}(\psi=1) \vee\max_{v \in\cB_0(k)}
\Pro^{
n}_v(\psi=0) &\geq& \Pro_{0}^{ n}(
\psi=1) \vee\max_{S \in
\mathcal{S}} \Pro^{ n}_{u(S)}(
\psi=0)
\\
&\geq& \Pro_{0}^{ n}(\psi=1) \vee\Pro^{ n}_{\mathcal{S}}(
\psi=0)
\\
&\geq& \frac{e^{1-e^{C_\nu}}}{4} \vee\frac{1-\sqrt{(e^{C_\nu
}-1)/2}}{2}=\frac{1}{2}-\nu,
\end{eqnarray*}
where the last inequality is a direct consequence of (\ref{EQchitot})
and Tsybakov (\citeyear{Tsy09}), Theorem 2.2, case (iii).
\end{pf*}

We observe a gap between our upper and lower bounds, with a term in
$\log(p/k)$ in the upper bound, and one in\vadjust{\goodbreak} $\log(p/k^2)$ in the lower
bound. This gap has been observed in the detection literature before
[see, e.g., \citet{Bar02,Ver10}, for an explicit remark] and, to our
knowledge, has never been addressed. However, if $p \geq
k^{2+\varepsilon}$, $\varepsilon>0$, upper and lower bounds match up
to constants, and the detection rate for the sparse eigenvalue is
optimal in a minimax sense. Under this assumption, detection becomes
impossible if $\theta< C\sqrt{(k/n) \log(p/k)}$
for a small enough constant $C>0$.

%s6 #&#
\section{Efficient methods for sparse principal component testing}
\label{SECSDP}

Computing the largest $k$-sparse eigenvalue $\lambda^k_{\max}$ of a
symmetric matrix $A$ is, in general, a hard computational problem. To
see this, consider the particular case where $A$ is a $p \times p$
symmetric matrix with values in $\{0,1\}$ and $A_{ii}=1$ for all
diagonal entries, so that $A$ corresponds to the adjacency matrix of an
undirected graph. It is not hard to see that $\lambda_{\max}^k(A)\le
k$, with equality if and only if the graph of $A$ contains a clique of
size $k$. It is a well-known fact of computational complexity [\citet
{Kar72}] that the decision problem associated to finding whether a
graph contains a clique of size $k$ is NP-complete.
%Note that if $k$ were fixed, this problem would actually be polynomial
%in the size $p$ of the graph since there are ``only'' ${p \choose k}
%in $k$ is clearly an issue even for moderate values of $k$.

%s6.1 #&#
\subsection{\texorpdfstring{Semidefinite relaxation for $\lambda^k_{\max}$}
{Semidefinite relaxation for lambda k max}}

Semidefinite programming (SDP) is the matrix equivalent of linear
programming. Define the scalar product in $\Sy_d$ by $\langle A,B
\rangle= \tr(AB)$. A semidefinite program can be written in the
canonical form:
%
%e6.1 #&#
%
\begin{eqnarray}\label{EQSDPcan}
\SDP&=& \mathrm{max.}\hspace*{31.4pt}\tr(CX)
\nonumber\\
&&\mbox{subject to}\quad \tr(A_i X) \leq b_i\qquad \forall i
\in\{1,\ldots, m\},
\\
&&\qquad\hspace*{33.05pt} X \succeq0.
\nonumber
\end{eqnarray}

As convex problems, they are computationally efficient and can be
solved using interior point or first order methods; see, for
example, \citet{BoyVan04,NesNem87}. Using SDP relaxations of problems
with nonconvex constraints such as integer programs is a common method
to find approximate solutions. Approximation bounds, up to a constant,
can sometimes be proved as in the celebrated result of \citet{GoeWil95}
for the MAXCUT problem. A major breakthrough for sparse PCA was
achieved by \citet{dAsGhaJor07}, who introduced a SDP relaxation for
$\lambda_{\max}^k$, but tightness of this relaxation is, to this day,
unknown. Our task is not as difficult though. Indeed, we only need to
prove that the SDP objective criterion has significantly different
behavior under $H_0$ and $H_1$.

Making the change of variables $Z=xx^\top$ in (\ref{EQlambdakmax2})
yields
\begin{eqnarray*}
\lambda^k_{\max}(A) &=& \mathrm{max.}\hspace*{31.4pt} \tr(AZ)
\\
&&\mbox{subject to}\quad \tr(Z) =1, |Z|_0 \leq k^2,
\\
&&\qquad\hspace*{33.05pt} Z \succeq0, \rk(Z)=1.
\end{eqnarray*}
Note that this problem contains two sources of nonconvexity: the $\ell
_0$ norm constraint and the rank constraint. We make two relaxations in
order to have a convex feasible set.
First, for a semidefinite matrix $Z$, with trace 1, and sparsity $k^2$,
the Cauchy--Schwarz inequality yields $|Z|_1 \leq k$, which is
substituted to the cardinality constraint in this relaxation.
Simply
dropping the rank constraint leads to the following relaxation of our
original problem:
%
%e6.2 #&#
%
\begin{eqnarray}\label{EQsdp}
\SDP_k(A)&=& \mathrm{max.}\hspace*{31.4pt} \tr(AZ)
\nonumber\\
&& \mbox{subject to}\quad \tr(Z) =1, |Z|_1 \leq k,
\\
&&\qquad\hspace*{33.05pt} Z \succeq0.
\nonumber
\end{eqnarray}
Note that this optimization problem is convex since it consists in
minimizing a linear objective over a convex set. Moreover, it is a
standard exercise to show that it can be expressed in the canonical
form (\ref{EQSDPcan}). As such, it can be solved efficiently using
any of the aforementioned algorithms. This natural relaxation was
originally developed in \citet{dAsGhaJor07}. Note that building on an
earlier version of this paper, \citet{dAsBacGha12} proposed a new SDP
relaxation to the same problem and derive somewhat larger detection
levels, at least for the interesting case where $k$ is small compared
to $p$.

Let us now study the behavior of the objective value $\SDP_k(\hat
\Sigma)$ under $H_1$ and~$H_0$, respectively. First, as a relaxation of
the original problem, for any $A \succeq0$, it holds
%
%e6.3 #&#
%
\begin{equation}
\label{EQlbSDP} \lambda^k_{\max}(A) \leq
\SDP_k(A).
\end{equation}
Since we have proved in Section~\ref{SECmain} that $\lambda^k_{\max
}(\Sh)$ takes large values under $H_1$, this inequality tells us that
so does $\SDP_k(\Sh)$. It remains to show that it stays small under
$H_0$. This can be achieved by using the dual formulation of the SDP.

%le6.1 #&#
%
\begin{lem}[{[\citet{BacAhidAs10}]}]
\label{LEMdualsdp}
For a given $A \succeq0$, we have by duality
\[
\SDP_k(A) = \min_{U \in\Sy_p} \bigl\{
\lambda_{\max}(A+U) + k|U|_{\infty} \bigr\}.
\]
\end{lem}

Together with (\ref{EQlbSDP}), Lemma~\ref{LEMdualsdp} implies that
for any $z \ge0$ and any matrix $U \in\bS_p$ such that $|U|_\infty
\le z$, it holds
%
%e6.4 #&#
%
\begin{equation}
\label{EQdualsdp} \lambda^k_{\max}(A) \leq
\SDP_k(A) \le\lambda_{\max}(A+U) + kz.
\end{equation}
A direct consequence of (\ref{EQdualsdp}) is that the functional
$\lambda^k_{\max}$ is robust to small perturbations in $|\cdot
|_\infty$-norm. Let $A \succeq0$ be such that its largest eigenvector
is $k$ sparse. Then, for any matrix $N$, (\ref{EQdualsdp}) yields
\[
\lambda^k_{\max}(A+N) \leq\lambda_{\max}
\bigl((A+N)-N\bigr) + k|N|_{\infty
} = \lambda_{\max}^k(A)+k|N|_{\infty}.
\]

%s6.2 #&#
\subsection{High probability bounds for convex relaxation}

We now study the properties of $\SDP_k(\hat{\Sigma})$ and other
computationally efficient variants as test statistics for our detection problem.
In view of (\ref{EQlbSDP}), the following proposition follows
directly from Proposition~\ref{THH1lk}.

%pr6.1 #&#
%
\begin{prop}
\label{THH1sdp}
Under $H_1$, we have, with probability $1-\delta$
\[
\SDP_k (\hat{\Sigma}) \geq1+ \theta- 2(1+ \theta)\sqrt{
\frac
{\log(1/\delta)}{n}}.
\]
\end{prop}
%
%Akin to Proposition~\ref{THH1lk}, Proposition~\ref{THH1sdp} shows
%that $H_1$ is the easy case. Indeed, under $H_1$, the lower deviations
%of $\SDP_k (\hat{\Sigma})$ remain small
%and do not depend on $k$ or $p$.
We now turn to the upper deviations under $H_0$.
%
%pr6.2 #&#
%
\begin{prop}
\label{THH0sdp}
Under $H_0$, we have, with probability $1-\delta$,
\begin{eqnarray*}
\SDP_k(\hat{\Sigma}) &\leq& 1+2 \sqrt{\frac{k^2\log(4p^2/\delta
)}{n}}+2
\frac{k \log(4p^2/\delta)}{n} \\
&&{}+ 2\sqrt{\frac{\log
(2p/\delta)}{n}}+ 2\frac{\log(2p/\delta)}{n}.
\end{eqnarray*}
\end{prop}
\begin{pf}
Let $st_z(A)$ be the soft-threshold of $A$, with threshold $z$, defined
by $(st_z(A))_{ij} = \operatorname{sign} (A_{ij}) (|A_{ij}|-z)_+$.
It follows from (\ref{EQdualsdp}) that for any $A \succeq0$,
%
%e6.5 #&#
%
\begin{equation}
\label{EQdualST} \SDP_k(A) \leq\lambda_{\max}
\bigl(st_z(A)\bigr) + kz.
\end{equation}
Let $\hat{\Delta} = \operatorname{diag}(\Sh)$ be the diagonal
matrix with the same diagonal entries as $\Sh$, and $\hat{\Psi} =
\Sh- \hat{\Delta}$ the matrix of its off-diagonal entries, so that
$\hat{\Sigma} = \hat{\Delta}+ \hat{\Psi}$. Since $\hat{\Psi}$
and $\hat{\Delta}$ have disjoint supports, it follows that
%
%e6.6 #&#
%
\begin{equation}
\label{EQdecompST} st_z(\hat{\Sigma}) = st_z(\hat{
\Delta}) + st_z(\hat{\Psi}).
\end{equation}

We first control the largest off-diagonal element of $\Sh$ by bounding
$|\hat{\Psi}|_{\infty}$ with high probability.
For every $i,j$, we have
\[
\hat{\Psi}_{ij}=\frac{1}{2} \Biggl[ \frac1n \sum
_{k=1}^n \biggl[\frac
{1}{2}(X_{ki}+X_{kj})^2
-1\biggr] - \frac1n \sum_{k=1}^n \biggl[
\frac
{1}{2}(X_{ki}-X_{kj})^2 -1\biggr]
\Biggr].
\]
Under $H_0$, we have $X \sim\N(0,I_p)$, so by Lemma~\ref{LEMLM}, it
holds for $t>0$ that
\[
\mathbf{P} \biggl(|\hat{\Psi}_{ij}| \geq2 \sqrt{\frac{t}{n}}+ 2
\frac{t}{n} \biggr) \leq4e^{-t}.
\]
Hence, by union bound on the off-diagonal terms, we get
\[
\mathbf{P} \biggl( \max_{i<j}|\hat\Psi_{ij}| \geq2
\sqrt{\frac
{t}{n}}+ 2\frac{t}{n} \biggr) \leq2p^2e^{- t}.
\]
Taking $t= \log(4p^2/\delta)$ yields that $|\hat{\Psi}|_{\infty
}\leq z$, with probability $1-\delta/2$, where
%
%e6.7 #&#
%
\begin{equation}
\label{EQoffdiag} z = 2 \sqrt{\frac{ \log(4p^2/ \delta)}{n}}+ 2\frac{
\log(4p^2/
\delta)}{n}.
\end{equation}

Note\vspace*{1pt} now that if we take $z$ as in (\ref{EQoffdiag}), then $st_z(\hat
{\Psi}) = 0$ on an event $\cE$ of probability $1-\delta/2$.
Furthermore, since $\hat{\Delta}$ is a nonnegative diagonal matrix,
then (\ref{EQdecompST}) yields that on the event $\cE$, it holds
%
%e6.8 #&#
%
\begin{equation}
\label{EQdiag2} \lambda_{\max}\bigl(st_z(\hat{\Sigma})
\bigr) = \lambda_{\max}\bigl(st_z(\hat{\Delta})\bigr) \le
\lambda_{\max}(\hat{\Delta}) = \max_{1\leq i \leq
p}\hat{
\Delta}_{ii}.
\end{equation}

Next, we control the largest diagonal element of $\Sh$ as follows.
We have by definition of $\hat{\Delta}$, for every $i=1,\ldots, p$
\[
\hat{\Delta}_{ii} = \frac{1}{n} \sum
_{j=1}^n X_{ji}^2.
\]
Applying Lemma~\ref{LEMLM} and a union bound over the $p$ diagonal
terms, we get
\[
\mathbf{P} \biggl(\max_{1 \leq i \leq p} \hat{\Delta}_{ii}
\geq1+2\sqrt{\frac{t}{n}} + 2\frac{t}{n} \biggr) \leq p
e^{-t}.
\]
Taking $t= \log(2p/\delta) $ yields with probability $1-\delta/2$,
%
%e6.9 #&#
%
\begin{equation}
\label{EQdiag} \max_{1 \leq i \leq p}\hat{\Delta}_{ii} \leq1+
2\sqrt{\frac{ \log
(2p/\delta)}{n}} + 2\frac{ \log(2p/\delta)}{n}.
\end{equation}

To conclude the proof of Proposition~\ref{THH0sdp}, observe
that (\ref{EQdualST}) implies that for all $z\ge0$, we have
\[
\SDP_k(\hat{\Sigma}) \leq\lambda_{\max}\bigl(st_z(
\hat{\Sigma})\bigr) + kz \le\lambda_{\max}\bigl(st_z(\hat{
\Delta})\bigr) + \lambda_{\max
}\bigl(st_z(\hat{\Psi})\bigr)
+ kz,
\]
where we used (\ref{EQdecompST}) and the triangle inequality for the
operator norm.

Putting together (\ref{EQdiag2}) and (\ref{EQdiag}) completes the
proof.
\end{pf}

%s6.3 #&#
\subsection{Hypothesis testing with convex methods}

Using the notation from Section~\ref{SECpb}, the results of the
previous subsection can be written as
\[
\Pro_{H_0}\bigl(\SDP_k(\hat{\Sigma})> \tilde{
\tau}_0\bigr) \leq\delta,\qquad \Pro_{H_1}\bigl(
\SDP_k(\hat{\Sigma})< \tilde{\tau}_1\bigr) \leq\delta,
\]
where $\tilde{\tau}_0$ and $\tilde{\tau}_1$ are given by
\begin{eqnarray*}
\tilde{\tau}_0 &=& 1+2 \sqrt{\frac{k^2\log(4p^2/\delta
)}{n}}+2\frac{k \log(4p^2/\delta)}{n}
+ 2\sqrt{\frac{\log
(2p/\delta)}{n}}+ 2\frac{\log(2p/\delta)}{n},
\\
\tilde{\tau}_1 &=& 1+ \theta- 2(1+ \theta)\sqrt{\frac{\log
(1/\delta)}{n}}.
\end{eqnarray*}

Whenever $\tilde\tau_1>\tilde\tau_0$, we take $\tau\in[\tilde
\tau_0,\tilde\tau_1]$ and define the following computationally
efficient test
$
\tilde\psi(\hat{\Sigma}) = \bone{\{ \SDP_k(\hat{\Sigma})> \tau
\}}.
$
It discriminates between $H_1$ and $H_0$ with probability $1-\delta$.\vadjust{\goodbreak}

It remains to find for which values of $\theta$ the condition $\tilde
\tau_1> \tilde\tau_0$ holds. It corresponds to our minimum detection level.

%th6.1 #&#
%
\begin{theorem}
\label{CORbarthetaSDP}
Assume that $p, n, k$ and $\delta$ are such that $\tilde{\theta}
\leq1$, where
%
%e6.10 #&#
%
\begin{eqnarray}
\label{EQbarthetaSDP} \tilde\theta&:=& 2 \sqrt{\frac{k^2\log(4p^2/\delta)}{n}}+2
\frac{k
\log(4p^2/\delta)}{n} +
2\sqrt{\frac{\log(2p/\delta)}{n}}\nonumber\\[-8pt]\\[-8pt]
&&{}+ 2\frac
{\log(2p/\delta)}{n}+ 4\sqrt{
\frac{\log({1}/{\delta})}{n}}.\nonumber
\end{eqnarray}
Then, for any $\theta> \tilde\theta$, any $\tau\in[\tilde\tau
_0,\tilde\tau_1]$, the test $\tilde\psi(\hat{\Sigma}) = \bone\{
\SDP_k(\hat{\Sigma})> \tau\}$ discriminates between $H_0$ and $H_1$
with probability $1-\delta$.
\end{theorem}
If we consider asymptotic regimes, for large $p,n,k$, taking $\delta=
p^{-\beta}$ with \mbox{$\beta>0$}, provides a sequence of tests $\tilde\psi
_n$ that discriminate between $H_0$ and $H_1$ with probability
converging to 1, for any fixed $\theta>0$, if $k^2 \log(p
)/n \rightarrow0$.

Note that,\vspace*{1pt} compared to Theorem~\ref{CORbartheta}, the price to pay
for using this convex relaxation is to multiply the minimum detection
level by a factor $\sqrt k$. Such a gap is observed for these
techniques in \citet{AmiWai09}. Nevertheless, in most examples, $k$
remains small and so is this price. As we will see in Section \ref
{SECCTLB}, there is strong evidence that $\tilde\tau_0$, which
dominates the detection rate, cannot be made smaller and that
therefore, our proof is tight.

%s6.4 #&#
\subsection{Simple methods}
\label{SUBsimple} While the SDP relaxation proposed in the previous
subsection is provably computationally efficient, it is also known to
scale poorly on large problems. Simple heuristics such as the diagonal
method of \citet{JohLu09} become more attractive for larger
problems. A~careful inspection of the proofs in the previous subsection
is quite informative. It indicates that our results not only hold for
the test $\tilde\psi(\Sh)$ but for a test based on a simpler statistic
arising from the dual formulation (\ref {EQdualST}). Indeed, to control
the behavior of $\SDP_k(\Sh)$ under~$H_0$, we showed that it was no
larger than the \emph{minimum dual perturbation} $\MDP_k(\Sh)$ defined
by
%
%e6.11 #&#
%
\begin{equation}
\label{EQmdp} \MDP_k(\Sh)=\min_{z\geq0} \bigl\{
\lambda_{\max}\bigl(st_z(\Sh)\bigr) + kz \bigr\}.
\end{equation}
Clearly $\MDP_k(\Sh) \ge\SDP_k(\Sh) \ge\lambda_{\max}^k( \Sh)$
so that both Propositions~\ref{THH1sdp}\break
and~\ref{THH0sdp} still hold for $\SDP_k(\Sh)$ replaced\vspace*{1pt} by $\MDP_k(\Sh)$.
As a result, for any $\theta> \tilde\theta$ the test $\hat\psi(\Sh
)=\bone\{ \MDP_k(\Sh) >\tau\}$ discriminates between $H_0$ and
$H_1$ with probability \mbox{$1-\delta$}.

Actually, a detection level of the same order as $\tilde\theta$ holds
already for an even simpler test statistic: the largest diagonal
element of $\Sh$. This method called \emph{Johnstone's diagonal
method} was first proposed by \citet{JohLu09} and later studied
by \citet{AmiWai09}. For the\vspace*{1pt} problem of detection considered here, it
dictates one to employ the test statistic
$\DD(\Sh) = \max_{1 \leq i \leq p} \Sh_{ii}$.
Using even simpler techniques than in Propositions~\ref{THH1sdp}
and~\ref{THH0sdp}, it is not hard to show that
\[
\Pro_{H_0}\bigl(\DD(\hat{\Sigma})> \tau^d_0
\bigr) \leq\delta,\qquad \Pro_{H_1}\bigl(\DD(\hat{\Sigma})<
\tau^d_1\bigr) \leq\delta
\]
for quantiles $\tau^d_0$ and $\tau^d_1$ given by
\begin{eqnarray*}
\tau^d_0 &=& 1 + \frac{1}{k} \theta- 2 \biggl(1+
\frac{1}{k} \theta\biggr)\sqrt{\frac{\log(1/\delta)}{n}},
\\
\tau^d_1 &=& 1 + 2 \sqrt{\frac{\log(p/\delta)}{n}}+2
\frac{ \log
(p/\delta)}{n}.
\end{eqnarray*}
However, as we shall see in Section~\ref{SECsim}, on simulated
datasets, $\MDP_k$ behaves much better than $\DD$ in practice. It was
proved by \citet{AmiWai09} that if the SDP (\ref{EQsdp}) has a
solution of rank one, then it is strictly better than Johnstone's
diagonal method. While they study a support recovery problem different
from the detection problem considered here, it seems to indicate that
the two methods are qualitatively different. However, the assumption
that the SDP (\ref{EQsdp}) has a solution of rank one is strong and
unnecessary in our problem. Indeed, our results from Section \ref
{SECCTLB} indicate that, if detecting a planted clique in a random
graph is computationally hard, then for large $(p,n,k)$, the SDP method
does not achieve better rates than the ones we prove. In particular,
this result is a good indication that with high probability, the
solution of the SDP is not rank-one for parameters in a range around
the minimax detection level.

%s7 #&#
\section{Generalization with weakened assumptions}
\label{SECweak}

In this section we investigate several extensions of our original
problem. For simplicity, we denote by $\starDP_k$ any of the two
functionals $\MDP_k$ or $\SDP_k$.

%s7.1 #&#
\subsection{\texorpdfstring{Sparsity in terms of $\ell_q$ norm}{Sparsity in terms of l q norm}}

Fix $q\in(0,2)$, and recall that $\cB_q(R)$ is the set of unit
vectors that are in an $\ell_q$ ball of radius $R>0$. This relaxed
notion of sparsity allows for vectors $v\in\R^p$ to have ordered
coordinates that decay fast enough but never take value zero. Note that
$q=2$ corresponds to no sparsity and requires different techniques. It
is therefore excluded from this section. Consider the following
hypothesis testing problem:
\begin{eqnarray*}
H_0 \dvtx X &\sim&\N(0,I_p),
\\
\tilde{H}_1^q \dvtx X &\sim&\N\bigl(0,I_p+
\theta vv^\top\bigr),\qquad v \in\cB_q\bigl(k^{{1}/{q} - {1}/{2}}
\bigr).
\end{eqnarray*}
The radius $k^{{1}/{q} - {1}/{2}}$ is the smallest $R>0$ such
that $\cB_0(k)\subset\cB_q(R)$, making it the most natural
relaxation of the notion of $k$-sparse vectors. Below, we show that it
yields the same detection levels as for $q=0$.
%
%th7.1 #&#
%
\begin{theorem}
\label{THinfminimaxlq}
Fix $\nu>0$. There exists a constant $C_\nu>0$ such that if
%
%e7.1 #&#
%
\begin{equation}
\label{EQinfminimaxlq} \theta< \underline{\theta}_\nu:=\sqrt{
\frac{ k\log(C_\nu
p/k^2+ 1 )}{n}} \wedge\frac{1}{\sqrt{2}},
\end{equation}
it holds, for $q\in(0,2)$
%
%e7.2 #&#
%
\begin{equation}
\inf_{\psi} \Bigl\{ \Pro_{0}^{ n}(\psi=1)
\vee\max_{v \in\cB
_q(k^{{1}/{q} -{1}/{2}})} \Pro^{ n}_v(\psi=0) \Bigr
\} \geq\frac12 - \nu,
\end{equation}
where the infimum is taken over all possible tests.
\end{theorem}
\begin{pf}
Let $v \in\R^p$ be a unit vector with sparsity $k$. It follows from
H\"older's inequality that $|v|_q \leq k^{1/q - 1/2}$.
%as a direct consequence of H\"older's inequality
%&\leq& \left( \sum_{j=1}^p |v_j|^2\right)^{\frac q2} k^{1-\frac q2}
Therefore, for any test $\psi$, we have
\[
\max_{v \in\cB_q(k^{{1}/{q} -{1}/{2}})} \Pro^{
n}_v(\psi=0) \geq
\max_{v \in\cB_0(k)} \Pro^{ n}_v(\psi=0),
\]
and the result follows as a direct consequence of Theorem~\ref{THinfminimax}.
\end{pf}
%
%As we will see using the radius $k^{\frac{1}{q} -\frac{1}{2}}$ give a
%lighter notation below. Nevertheless, solving the equation $\omega=k^{
%leads to the following lower bound on the optimal detection level when
%$v \in\cB_q(\omega)$ for any $\omega> 1$:

To show a matching upper bound, we use the following lemma.
%
%le7.1 #&#
%
\begin{lem}
\label{LEMsparselq}
Let $v \in\R^p$ be a unit vector, $|v|_2 =1$. Then, for any $r \geq
1$, there exists a $r$-sparse unit vector $x \in\cB_0(r)$ such that
\[
1- r^{1-{2}/{q}} |v|_q^2 \le\bigl(x^\top v
\bigr)^2 \le1.
\]
\end{lem}
\begin{pf}
Assume\vspace*{1pt} without loss of generality that $|v_1| \geq\cdots\geq|v_p|$.
Define $\tilde x_j=v_j$ if $j\le r$, $\tilde x_j=0$ otherwise, and
$x=\tilde x/|\tilde x|_2$.
We have $(x^\top v)^2 = \tilde x ^\top v = \sum_{j=1}^r |v_j|^2= 1-
\sum_{j=r+1}^p|v_j|^2$. Moreover, since $|v_r| \leq r^{- 1/q} |v|_q$,
\begin{eqnarray*}
\sum_{j=r+1}^p|v_j|^2
&\le&\sum_{j=r+1}^p |v_r|^{2-q}
|v_j|^q\le|v|_q^{2-q}
r^{1 -2/q} \sum_{j=r+1}^p
|v_j|^q\\
&\le& r^{1-
{2}/{q}} |v|_q^2.
\end{eqnarray*}
\upqed\end{pf}
Vectors in $\cB_q(k^{{1}/{q} - {1}/{2}})$ can therefore be
approximated by sparse unit vectors. This property can be leveraged to
show that for an appropriate choice of $k_q$, a test based on $\lambda
^{k_q}_{\max}(\hat\Sigma)$ is optimal.

%pr7.1 #&#
%
\begin{prop}
\label{THtH1lk}
Under $\tilde{H}^q_1$, let $\varepsilon\in(0,1)$, and define $k_q$
to be the smallest integer such that $k_q \ge k \varepsilon^{
{1}/({1-2/q})}$. Then with probability $1-\delta$,
\[
\lambda^{k_q}_{\max}(\hat{\Sigma}) \geq1+ (1-\varepsilon)
\theta- 2 (1+\theta) \sqrt{\frac{\log(2/\delta)}{n}}.
\]
\end{prop}
\begin{pf}
Let $x \in\R^p$ be the $k_q$-sparse unit norm approximation of $v$
from Lemma~\ref{LEMsparselq}. It follows from the proof of
Proposition~\ref{THH1lk} that
\[
\lambda^{k_q}_{\max}(\hat{\Sigma}) \geq1+ \theta
\bigl(v^\top x\bigr)^2 - 2 \bigl(1+\theta
\bigl(v^\top x\bigr)^2\bigr) \sqrt{\frac{\log(2/\delta)}{n}}.
\]
Lemma~\ref{LEMsparselq} with $r=k_q \ge k \varepsilon^{
{1}/({1-2/q})}$ yields $1-\eps\le(x^\top v)^2 \le1$.
\end{pf}
Moreover, it follows from Proposition~\ref{THH0lk} that for any
$\varepsilon\in(0,1)$ and integer $k_q$, with probability $1-\delta
$, it holds
\[
\lambda_{\max}^{k_q}(\hat{\Sigma}) \leq1+ 2 \biggl(
\frac{k_q
\log(9ep / k_q)+\log(2/\delta)}{n} + \sqrt{\frac{k_q \log(9ep /
k_q)+\log(2/\delta)}{n}} \biggr).
\]

Since $k_q$ is only a constant factor away from $k$ for all $q \in
(0,2)$ and $\varepsilon\in(0,1)$, the statistic $\lambda_{\max
}^{k_q}(\hat{\Sigma})$ achieves optimal rates of detection.

In an estimation context, \citet{VuLei12} [see also \citet
{PauJoh12,BirJohNadPau12} for related results using a different method]
have examined the $\ell_q$ sparsity assumption for $q \in(0,1]$.
Their estimation method consists of maximizing the quadratic form $x
\mapsto x^\top\hat\Sigma x$ over $\cB_q(R)$ for some given $R>0$. We
argue that in light of Lemma~\ref{LEMsparselq}, the estimation
problem of \citet{VuLei12} can be solved by maximizing the quadratic
form over $\cB_0(R')$ for some appropriate choice of $R'$ that depends
on $k$ and $q$ and extended to $q \in[0,2)$. In particular, an
algorithm for $\ell_0$-sparse PCA can be used for $\ell_q$-sparse PCA.

Similar results hold for our convex relaxations. Following the same
steps as in the proof of Theorem~\ref{CORbarthetaSDP}, we find that
there exists a constant $C_{q}>0$ such that tests based on $\starDP
_{k_q}$ discriminate between $H_0$ and $\tilde{H}^q_1$ with
probability $1-\delta$, for any $\theta> C_{q} \tilde\theta$, where
$\tilde\theta$ is defined in (\ref{EQbarthetaSDP}). In particular,
a gap of size $\sqrt{k}$ is observed between these methods and the
optimal ones.
%Nevertheless, for $q = 1$, since the SDP relaxation employs an $
%if $v \in\cB_1(k^{1/2})$, the matrix $Z=vv^\top$ satisfies the
%constraint of \eqref{EQsdp}. Therefore, under $\tilde H^q_1$, it
%holds with probability $1-\delta$ that
%$$
%$$
%Together with Proposition~\ref{THH0sdp} and its extension in
%subsection~\ref{SUBsimple}, it implies that the tests based on $

%Assume that $p, n, k,$ and $\delta$ are such that $\tilde\theta\le
%1$.
%Then, for any $\theta> \tilde\theta$ and for any $\tau\in[\tilde
%probability $1-\delta$.

%s7.2 #&#
\subsection{Sub-Gaussian random variables}

Our results can be extended to more general assumptions, where the
variables $X_1,\ldots,X_n \in\R^p$ are sub-Gaussian in the following sense.

%de7.1 #&#
%
\begin{defin}
A real-valued random variable $G$ is said to be standard sub-Gaussian if
$
\E[\exp(t(G - \E[G]) ) ] \le\exp(t^2/2)$ for
all $t\in\R$.
\end{defin}

Let $Z_1,\ldots,Z_n \in\R^p$ be i.i.d. vectors with i.i.d. standard
sub-Gaussian coefficients, such that for all $i = 1,\ldots, n$ it holds
$\E[Z_i] = 0, \E[Z_i Z_i^\top] = I_p$.

Given a \emph{scatter matrix} $\Sigma\succeq0$, for any $i=1,\ldots,
n$, define $X_i = \Sigma^{1/2} Z_i$. Sub-Gaussian random vectors were
generated in the same way by \citet{VuLei12}. Under this condition, we
define the new detection problem with hypotheses $H'_0$ and $H'_1$, for
$\theta>0$ by
\begin{eqnarray*}
H'_0 \dvtx\Sigma&=& I_p,
\\
H'_1 \dvtx\Sigma&=& I_p+ \theta v
v^\top,\qquad v \in\cB_0(k).
\end{eqnarray*}
Replacing Lemma~\ref{LEMLM} by Lemma~\ref{LEMLMsubexp} in the
proofs of Propositions~\ref{THH0lk} and~\ref{THH0sdp},
we get, respectively, the two following results.

%pr7.2 #&#
%
\begin{prop}
Under $H'_1$, for $\theta\le1$, it holds with probability $1-\delta$
\[
\lambda^k_{\max}(\hat{\Sigma}) \geq1+ \theta- 6 \biggl(64
\frac{
\log(2/\delta)}{n} +32\sqrt{\frac{\log(2/\delta)}{n}} \biggr).
\]
Moreover, under $H'_0$, it holds with probability $1-\delta$,
\[
\lambda_{\max}^k(\hat{\Sigma}) \leq1+ 352 \biggl( 2
\frac{k \log
(9ep / k)+\log(2/\delta)}{n} + \sqrt{\frac{k \log(9ep / k)+\log
(2/\delta)}{n}} \biggr).
\]
\end{prop}
Similarly, for the $\starDP_k$ statistic, we obtain the following bound.

%pr7.3 #&#
%
\begin{prop}
\label{THH0sdpsg}
Under $H_0$, we have, with probability $1-\delta$,
\[
\SDP_k(\hat{\Sigma})\le\MDP_k(\hat{\Sigma}) \leq1+6
\biggl( 64\sqrt{\frac{k^2\log(4p^2/\delta)}{n}}+128\frac{k \log
(4p^2/\delta)}{n} \biggr).
\]
\end{prop}
As a consequence, all the results from Sections~\ref{SECmain}
and~\ref{SECSDP} can be extended to the present sub-Gaussian case.
In particular, the same gap between the detection levels of the two
procedures is observed.

%s7.3 #&#
\subsection{Adversarial noise}

While our previous results rely heavily on the fact that the $X_i$ are
sub-Gaussian random vectors, we can find much weaker assumptions under
which the results for detection using the $\starDP$ statistics are
still valid. We also describe an adversarial noise setting in which the
detection level attained by $\starDP_k$ is actually optimal. Assume that
%
%e7.3 #&#
%
\begin{equation}
\label{EQadv} \Sh= \Sigma+ N.
\end{equation}
Here the only assumption on $N$ is that $|N|_{\infty} \leq\sqrt{\log
(p/\delta)/n}$ with probability $1-\delta$. Up to constant factor,
this is a generalization of our initial setting, and can describe a
situation where the data is censured, akin to the setting of \citet
{LohWai12}, but where the censured entries are not necessarily chosen
at random.

%We show below that the high probability bounds previously obtained for
%$ \starDP_k(\Sh)$ depend only on this mild assumption.
%
%pr7.4 #&#
%
\begin{prop}
\label{THH1adv}
Under $H_1$, we have with probability $1-\delta$
\[
\starDP_k(\Sh) \geq\lambda^k_{\max}(\hat{
\Sigma}) \geq1+ \theta- k\sqrt{\frac{\log(p/\delta)}{n}}.\vadjust{\goodbreak}
\]
\end{prop}
\begin{pf}
Recall that for any $v$ such that $|v|_0 \leq k$, we have
\begin{eqnarray*}
\starDP_k(\Sh) &\geq&\lambda^k_{\max}(\hat{
\Sigma}) \geq v^\top\Sh v\geq v^\top\bigl(I_p+
\theta vv^\top\bigr)v +v^\top Nv
\\
&\geq& 1 + \theta-|N|_{\infty}|v|_1^2\ge1 +
\theta- k |N|_{\infty
},
\end{eqnarray*}
which yields the desired result.
\end{pf}
%
%pr7.5 #&#
%
\begin{prop}
\label{THH0adv}
Under $H_0$, we have with probability $1-\delta$
\[
\lambda^k_{\max}(\hat{\Sigma}) \leq\starDP_k(
\Sh) \leq1+ k\sqrt{\frac{\log(p/\delta)}{n}}.
\]
\end{prop}
\begin{pf}
It follows from (\ref{EQdualsdp}) that
$
\lambda^k_{\max}(\hat{\Sigma}) \leq\starDP_k(\Sh) \leq\lambda
_{\max}(I_p) + k|N|_{\infty},
$
which yields the desired result.
\end{pf}

The following theorem follows from Propositions~\ref{THH1adv} and
\ref{THH0adv}. We omit its proof.
%
%th7.2 #&#
%
\begin{theorem}
\label{CORtestadv}
Let $\psi^{\mathrm{adv}}$ be the test defined by
\[
\psi^{\mathrm{adv}}(\Sh) = \bone{ \biggl\{\starDP_k(\Sh) > 1 + k
\sqrt{\frac
{\log(p/\delta)}{n}} \biggr\}}.
\]
Then $\psi^{\mathrm{adv}}$ discriminates between $H_0$ and $H_1$ with
probability $1-\delta$ if
$
\theta> 2k\sqrt{\log(p/\delta)/n}.
$
\end{theorem}

We now prove the corresponding lower bound. Let $v = (v_1,\ldots,
v_p)^\top\in\R^p$ be such that $v_j=1/\sqrt k$ if $j \le k$ and
$v_j=0$ otherwise. Define the random matrix $N$ that takes values $\pm
\frac{\theta}{2}vv^\top$, each with probability $1/2$.

%th7.3 #&#
%
\begin{theorem}
\label{THinfadv}
There exists an adversarial model of the form (\ref{EQadv}) where
$|N|_\infty\le\sqrt{\log(p)/n}$, such that if $\theta\le2k\sqrt{\log
(p)/n}$, then for any test $\psi(\Sh) \in\{0,1\}$ it holds
\[
\Pro_{H_1}\bigl(\psi(\Sh)=0\bigr) \vee\Pro_{H_0}\bigl(\psi(
\Sh)=1\bigr) \ge1/2.
\]
\end{theorem}
\begin{pf}
Note first that $|N|_{\infty} = \theta/(2k) \leq\sqrt{(\log p)/n}$
so that
\[
\Pro_{H_0} \biggl(\Sh= I_p + \frac{\theta}{2}
vv^\top\biggr) = \frac
{1}{2},\qquad \Pro_{H_1} \biggl(\Sh=
I_p + \frac{\theta}{2} vv^\top\biggr) =
\frac{1}{2}.
\]
Therefore, if $\psi( I_p + \frac{\theta}{2} vv^\top)=1$, then $\Pro
_{H_0}(\psi(\Sh)=1) \ge1/2$ and if $\psi( I_p + \frac{\theta}{2}
vv^\top)=0$, then $\Pro_{H_1}(\psi(\Sh)=0) \ge1/2$.
\end{pf}
Note that the lower bound in Theorem~\ref{THinfadv} below is not
minimax since there exists one model under which all tests cannot
discriminate between $H_0$ and $H_1$ with probability less than $1/2$.
It implies that tests based on either $\starDP_k$ and $\lambda_{\max
}^k$ are optimal.

%s8 #&#
\section{Complexity theoretic lower bounds}
\label{SECCTLB}
The difference between the detection rates proved for the testing
statistic $\lambda_{\max}^k$ and the convex optimization based
statistics $\SDP
_k$ and $\MDP_k$ suggests a statistical cost for computational
efficiency. Such phenomena are hinted at by \citet{ChaJor13}. While it
is not hard to see that our bounds are tight for the diagonal method,
it is legitimate to wonder if the observed gap for $\SDP_k$ and $\MDP
_k$ comes from a proof artifact, or an intrinsic limitation of the
problem. The computational hardness of the related RIP certification
has recently attracted a lot of interest. By reductions to problems
with known complexity theoretic limitations, \citet{BanDob12} and \citet
{KoiZou12} prove that it is in general impossible to approximate in
polynomial time the $\lambda_{\max}^k$ statistic up to an arbitrarily small
constant. Clearly, a constant factor approximation of $\lambda_{\max
}^k$ would
suffice to achieve optimal detection rates. However, such results are
not sufficient for two reasons. First they do not rule out the
existence of a polynomial time algorithm that approximates $\lambda
_{\max}^k$
within a large enough constant. Second, such results are in nature
\emph{worst case}, meaning that the input matrix can be arbitrarily
difficult for an algorithm. Rather, in our problem, the entry matrix is
an empirical covariance matrix constructed from i.i.d. random vectors
with Gaussian distribution.
Hereafter, we develop a \emph{polynomial time reduction} from another
problem which is believed to be hard in average: the planted clique problem.

%s8.1 #&#
\subsection{Reduction to the planted clique problem}

A careful inspection of the proof of Corollary~\ref{CORbarthetaSDP}
and the results of Section~\ref{SUBsimple} reveals that the only
way to obtain better detection levels for the $\SDP$ and $\MDP$
statistics is to prove a better control of the statistics under the
null hypothesis. We argue below that this is unlikely.

Let $X_1,\ldots,X_n \in\R^p$ be i.i.d. Gaussian vectors with
distribution $\cN(0,I_p)$ and for any $\alpha\in[1,2]$, consider the
following hypothetical bound:
%
%e8.1 #&#
%
\renewcommand{\theequation}{\mbox{$B_\alpha$}}
\begin{equation}
\label{EQBalpha}
\SDP_k(\Sh) \leq1 +
C_\alpha\sqrt{\frac{k^\alpha\log(p/\delta )}{n}}\qquad\mbox{with
probability } 1- \delta,
\end{equation}
where $C_\alpha>0$.
Our prior results hinge on proving that $B_2$ holds. However, to
achieve minimax optimal detection rates, one would need to prove $B_1$.
Reasoning by contradiction, we examine the consequences of $B_\alpha$
with $\alpha\in(1,2)$. In particular, such bounds would yield
polynomial time algorithms to detect small planted cliques in random
graphs. Hereafter, we argue that the existence of such algorithms is unlikely.

Fix an integer $k \ge0$, and let $\cG(n,1/2,k)$ be the distribution
over the set of graphs on $n$ vertices generated as follows. Pick $k$
vertices at random, and place a clique\footnote{A clique is a subset
of fully connected vertices.} between them; then connect every other
pair of vertices by an edge independently with probability\vadjust{\goodbreak} $1/2$. Note
that for $k=0$, $\cG(n,1/2,0)=\cG(n,1/2)$ is simply the distribution
of an Erd{\H{o}}s--R{\'e}nyi random graph. In the decision version of
the planted clique problem, called \textsf{Planted Clique}, one is
given a
graph $G$ on $n$ vertices and the goal is to test
\begin{eqnarray*}
H_0^{\mathsf{PC}} \dvtx G &\sim&\cG(n,1/2),
\\
H_1^{\mathsf{PC}} \dvtx G &\sim&\cG(n,1/2,k)
\end{eqnarray*}
for some given $k \ge2$ with probability of error at most $\delta>0$.
The search version of this problem consists of finding the clique
planted under $H_1$. The search problem was introduced by \citet{Jer92}
and \citet{Kuc95} while the decision version is traditionally
attributed to Saks; see \citet{KriVu02}, Section 5. It is known [see,
e.g., \citet{Spe94}] that if $k > 2\log_2 n$, the planted clique under
$H_1$ is the only clique of size $k$ in the graph, asymptotically
almost surely. We consider only such values of $k$ hereafter.

For $k=o(\sqrt{n})$ there is no known polynomial time algorithm that
solves this problem. The first polynomial time algorithm for the case
$k=C\sqrt{n}$ was proposed in \citet{AloKriSud98} and is based on
spectral techniques. Subsequent algorithms with similar performance
appeared in \citet{AmeVav11,DekGurPer11,FeiRon10,FeiKra00}. It is widely
believed that there is no polynomial time algorithm that solves
\textsf{Planted Clique} for any $k$ of order $n^c$ for some fixed
positive $c<1/2$, and it can even be proved that certain algorithmic
techniques such as the Metropolis process [\citet{Jer92}] and the
Lov\`asz--Schrijver hierarchy of relaxations [\citet{FeiKra03}] fail at
this task. Moreover, \textsf{Planted Clique} is provably hard in
certain computational models, as seen in \citet{Ros10},
\citet{FelGriRey13} which brings more evidence toward its hardness. Note
that recent results by \citet{BruVem09,FriKan08} based on
$r$-dimensional tensors, suggest an algorithmic approach capable of
finding a planted clique of size $O(n^{1/r})$, but currently this
tensor-based approach is not known to yield a polynomial time algorithm
for $r >2$. The confidence in the difficulty of this problem has led
researchers to prove hardness results assuming that the planted clique
problem is indeed hard. Examples include cryptographic applications
[\citet{JuePei00}], testing for $k$-wise dependence
[\citet{AloAndKau07}], approximating Nash equilibria [\citet{HazKra11}]
and approximating solutions to the densest $k$-subgraph problem
[\citet{AloAroMan11}].

Consider the following polynomial-time reduction from a graph instance
to random vectors, valid for the case $p=n$. Let $A$ be the $n \times
n$ adjacency matrix of a random graph $G$, and let $U$ be the $n \times
n$ matrix defined for any $1\le i \le j$ by
\[
U_{ij} = \cases{ 2A_{ij}-1, &\quad if $i<j$,
\cr
\varepsilon_{ij}, &\quad otherwise,}
\]
where $\{\eps_{ij}\}_{i,j}$ is a sequence of i.i.d. Rademacher $\pm1$
random variables. Moreover, let $Z^{(1)},\ldots, Z^{(n)} \in\R^n$
be $n$ i.i.d. $\cN(0, I_n)$ random vectors, and define $X_{ij} =
|Z_{j}^{(i)}|U_{ij}$. Finally define the $n \times n$ empirical
covariance matrix $\hat\Sigma$ associated to the vectors
$X_i=(X_{i1},\ldots, X_{in})^\top\in\R^n$ as in (\ref
{EQdefEmpCov}). This construction clearly takes polynomial time.

If $G \sim\cG(n,1/2)$, by construction, $X_1,\ldots,X_n \in\R^n$
are i.i.d. centered standard Gaussian vectors, where all the
coefficients are independent. If $G \sim\cG(n,1/2,k)$, it is no
longer the case, and the $\starDP_k$ statistic behaves in a
qualitatively different manner.

%s8.2 #&#
\subsection{Computational theoretic lower bounds for $\SDP$ and $\MDP$}
In this subsection, we illustrate the intrinsic limitations of the
$\SDP$ and $\MDP$ methods in the detection problem using arguments
borrowed from computational complexity. We begin by showing that both
statistics take large values on the problem reduced from a graph with a
planted clique.
%
%le8.1 #&#
%
\begin{lem}
\label{LEMLBH1reduc}
Let $G \sim\cG(n,1/2,k)$, $k \ge14$ even, and $X_1,\ldots,X_n \in
\R^n$ be constructed as above. It holds, with probability $1-\delta$,
\[
\MDP_k (\Sh)\ge\SDP_k (\Sh) \ge1+ \frac{k^2}{4\pi n} - 3
\sqrt{\frac{k\log(2/\delta)}{n}}.
\]
\end{lem}
\begin{pf}
Let $S \subset\{1,\ldots, n\}$ be the random subset of $k$ vertices
on which the clique has been planted. By construction, there are
subsets $S_1$ and $S_2$ of $S$, of cardinality $k/2$, such that the
random variables $X_{ij}, i \in S_1, j \in S_2$ are all positive almost
surely. Assume without loss of generality that $S=\{1,\ldots,k\}$,
$S_1=\{ k/2 +1,\ldots,k\}$ and $S_2=\{1,\ldots, k/2 \}$. Let\vspace*{1pt}
$v=v(S_2)$ be the unit vector with $j$th coordinate equal to $2/\sqrt
{k}$ if $j \in S_2$ and $0$ otherwise. It yields
%
%e8.2 #&#
%
\setcounter{equation}{0}
\renewcommand{\theequation}{\arabic{section}.\arabic{equation}}
\begin{eqnarray}\label{EQSDPLB}
\MDP_k(\Sh)&\ge&\SDP_k(\Sh) \geq v^\top\Sh v =
\frac{1}{n}\sum_{i=1}^n
\bigl(v^\top X_i\bigr)^2
\nonumber
\\
&\geq&\frac1n \sum_{i\in S_1} \bigl(v^\top
X_i\bigr)^2 +\frac1n \sum_{i\notin
S}
\bigl(v^\top X_i\bigr)^2
\nonumber\\
&=&\frac1n \sum_{i\in S_1} \frac{1}{|S_2|} \biggl(\sum
_{j \in S_2} \bigl|Z_j^{(i)} \bigr|
\biggr)^2 +\frac{n-k}{n} \frac{1}{n-k} \sum
_{i\notin S} \bigl(v^\top X_i
\bigr)^2.
\end{eqnarray}
We begin by controlling the first term on the right-hand side of (\ref
{EQSDPLB}). For all $i \in S_1$, define the centered sub-Gaussian
random variable
\[
Y_i = \sum_{j \in S_2} \bigl\{
\bigl|Z_j^{(i)} \bigr|- \sqrt{2/\pi} \bigr\}
\]
and observe that
\begin{eqnarray*}
\frac1n \sum_{i\in S_1} \frac{1}{|S_2|} \biggl(\sum
_{j \in S_2} \bigl|Z_j^{(i)} \bigr|
\biggr)^2 &=& \frac1n \sum_{i\in S_1}
\frac{1}{|S_2|} \bigl(Y_i + |S_2| \sqrt{2/\pi}
\bigr)^2
\\
&\ge& \frac2 \pi\frac{|S_1|\cdot |S_2|}{n} + 2 \sqrt{\frac2 \pi} \frac1n
\sum
_{i\in S_1} Y_i.
\end{eqnarray*}
It follows from Lemma~\ref{LEMabsnorm} that with probability
$1-\delta/2$, we have
\[
\sum_{i\in S_1}Y_i \ge-
\sqrt{2|S_1|\cdot|S_2|\log(2/\delta)}.
\]
Together, the previous two displays yield
%
%e8.3 #&#
%
\begin{equation}
\label{EQpr91} \frac1n \sum_{i\in S_1}
\frac{1}{|S_2|} \biggl(\sum_{j \in S_2}
\bigl|Z_j^{(i)} \bigr| \biggr)^2 \ge\frac{k^2}{2\pi n}-
\frac{2k}{\sqrt{\pi
}n}\sqrt{\log(2/\delta)}.
\end{equation}

To control the second term on the right-hand side of (\ref
{EQSDPLB}), we use Lem\-ma~\ref{LEMLM}. It holds with probability
$1-\delta/2$ that
%
%e8.4 #&#
%
\begin{equation}
\label{EQpr92}\quad \frac{1}{n-k} \sum_{i\notin S}
\bigl(v^\top X_i\bigr)^2 = 1+
\frac{1}{n-k} \sum_{i\notin S} \bigl[
\bigl(v^\top X_i\bigr)^2 -1 \bigr] \ge1 - 2
\sqrt{\frac{\log(2/\delta)}{n-k}}.
\end{equation}
Therefore, with probability $1- \delta$, we get from (\ref
{EQSDPLB}), (\ref{EQpr91}) and (\ref{EQpr92}) that
\begin{eqnarray*}
\SDP_k(\Sh) & \ge& \frac{k^2}{2\pi n}-\frac{2k}{\sqrt{\pi}n}\sqrt{\log(2/
\delta)} +\frac{n-k}{n} \biggl( 1 - 2 \sqrt{\frac{\log
(2/\delta)}{n-k}} \biggr)
\\
&\ge& 1+ \frac{k^2}{4\pi n} - 3\sqrt{\frac{k\log(2/\delta)}{n}},
\end{eqnarray*}
where the last inequality holds for $k \ge14$.
This yields the desired result.
\end{pf}

Next, we prove that improving substantially the bound of
Proposition~\ref{THH0sdp} (i.e., if $B_\alpha$ were to hold for some
$\alpha\in[1,2)$) would allow us to detect the presence of cliques of
size $n^{c}$ for some $c<1/2$.
%
%th8.1 #&#
%
\begin{theorem}
\label{THcomp}
Let $X_1,\ldots,X_n \in\R^n$ be i.i.d. $\cN(0,I_n)$ random vector
and let $\Sh$ be their corresponding empirical covariance matrix as
defined in (\ref{EQdefEmpCov}). If for any $\alpha\in[1,2]$,
$(B_\alpha)$ is valid for $p=n$, that is,
%
%e8.5 #&#
%
\renewcommand{\theequation}{\mbox{$B_\alpha$}}
\begin{equation}
\label{EQUBSDPcomp}
\SDP_k(\Sh)\leq1 + C_\alpha\sqrt{\frac{k^\alpha\log(n/\delta
)}{n}}\qquad
\mbox{with probability } 1-\delta,
\end{equation}
where $C_\alpha>0$, then there exists a polynomial time algorithm that
discriminates between $\cG(n,1/2)$ and $\cG(n,1/2,k)$ with
probability $1-\delta$, as soon as $
k \ge[C n \log(n/\delta)]^{{1}/({4-\alpha})}$
for some constant $C>0$ that depends only on $C_\alpha$. The same
holds if $\SDP_k (\Sh)$ is replaced by $\MDP_k (\Sh)$. In
particular, for any fixed $\alpha<2$ and $\delta>0$, it allows one to
detect the presence of cliques of size $n^{c}$ for some $c<1/2$ with
probability $1-\delta$.
\end{theorem}
\begin{pf}
Note first that since $\SDP_k (\Sh) \le\MDP_k (\Sh)$, it suffices
to prove the result for $\SDP_k (\Sh)$.

Let $G$ be a random graph from $H_0^{\mathsf{PC}}$ or $H_1^{\mathsf
{PC}}$. Our
goal is to construct a test $\phi$ that discriminates between the two
hypotheses. Let $X_1,\ldots, X_n \in\R^n$ be $n$ random vectors
obtained by the polynomial time reduction described in the previous
subsection, and denote by $\Sh$ their associated empirical covariance
matrix. We propose the following test:% $\phi=0$ if
%Formally
%
\[
\phi=\phi(\Sh) = \mathbf{1} \biggl\{\SDP^{(\varepsilon)}_k(\Sh) > 1+
C_\alpha\sqrt{\frac{k^\alpha\log(n/\delta)}{n}} \biggr\},
\]
where $\SDP^{(\varepsilon)}_k(\Sh) \ge\SDP_k(\Sh) -\eps$ is an
approximation of the SDP solution with tolerance $\varepsilon\le
1/\sqrt{n}$. In particular, $\SDP^{(\varepsilon)}_k(\Sh)$ and thus
$\phi$ can be computed in polynomial time.

Recall that under $H_0$ (no planted clique), the $X_i$'s are i.i.d. $\cN
(0, I_n)$ so that \mbox{$\phi=0$} with probability $1-\delta$, which
controls the type I error appropriately. Moreover, we know from
Lemma~\ref{LEMLBH1reduc} that under $H_1$, we have $\phi=1$ with
probability $1-\delta$ as soon as
\[
1+ C_\alpha\sqrt{\frac{k^\alpha\log(n/\delta)}{n}} \le1+ \frac
{k^2}{4\pi n} - 3 \sqrt{
\frac{k \log(2/\delta)}{n}}.
\]
Solving for $k$ yields that it is sufficient to have $k \ge[Cn
\log(n/\delta) ]^{{1}/({4-\alpha})}$,
for some constant $C>0$ that depends only on $C_\alpha$. As a result,
our test allows us to detect the presence of cliques of size $ [Cn
\log(n/\delta) ]^{{1}/({4-\alpha})}$.
\end{pf}

The consequences of Theorem~\ref{THcomp} can be taken two ways. If
one believes that detecting planted cliques of size at most $O(n^c),
c<1/2$ is hard, then suboptimality by a factor $\sqrt{k}$ is intrinsic
to the SDP relaxation. Otherwise, the $\SDP_k$ statistic allows to
reach new detection levels for \textsf{Planted Clique}.

To conclude, observe that the above results apply to the specific tests
based on $\MDP$ and $\SDP$ only. An interesting question is to find
whether this limitation is intrinsic to \emph{all} polynomial time
computable tests. Currently, the main limitation of the above proof is
that $\SDP_k (\Sh)$ is well controlled under $H_1$, but it may no
longer be the case for any other statistic.

%s9 #&#
\section{Numerical experiments}
\label{SECsim}

Computational cost is a crucial element in this study. In \citet
{BacAhidAs10}, the SDP relaxation with accuracy $\varepsilon$ is shown
to have a total complexity of $\mathcal{O}(k p^3\sqrt{\log
(p)}/\varepsilon)$. This is achieved by minimizing a smooth\vadjust{\goodbreak}
approximation of the dual function, using first order methods from
\citet{Nes03}.
However, this polynomial cost is already prohibitive in a
high-dimensional setting, and we study only tests based on the $\MDP
_k$ statistic. The latter is the solution of a one-dimensional
minimization problem, and is approximately solved by taking a uniform
grid on the variable $z$.
The purpose of this section is to illustrate the empirical behavior of
tests based on $\MDP_k$ and to compare it with the diagonal method.

%s9.1 #&#
\subsection{Comparison of simple methods}

We simulate $N = 1000$ samples of $n$ independent random vectors
$X^0_1,\ldots,X^0_n\sim\mathcal{N}(0,I_p)$ and $X^1_1,\ldots,X^1_n
\sim\break\mathcal{N}(0,I_p + \theta vv^\top)$, for random unit vectors
$v$ supported on $S=\{1,\ldots, k\}$. The vector $v_S$ is distributed
uniformly on the unit sphere of dimension $k$.

%f1 #&#
%
\begin{figure}

\includegraphics{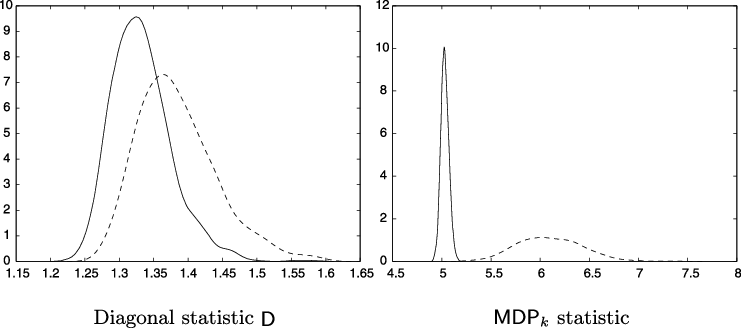}

\caption{For $p=500$, $n=200$, $k=30$, $N=1000$, estimated densities
for the two statistics, under~$H_0$ (whole line) and under $H_1$
(dashed line).}
\label{FIGkde}
\end{figure}

It yields $N$ empirical covariance matrices $\Sh_1^0,\ldots,\Sh_N^0$
under $H_0$ and $N$ of them, $\Sh_1^1,\ldots,\Sh_N^1$ under $H_1$. We
compute the $\DD$ and $\MDP_k$ statistics for these samples and
compare their densities. We take $\theta=4$ and observe that the $\DD
$ statistic yields two distributions under $H_0$ and $H_1$ that are
hard to distinguish (Figure~\ref{FIGkde}, left). In particular, it is
clear that the statistic $\DD$ cannot discriminate between $H_0$ and
$H_1$ for $\theta=4$, with this set of parameters. However, the
distributions of $\MDP_k(\Sh)$ under $H_0$ and $H_1$ have almost
disjoint support so that it can discriminate between the two hypotheses
with probability close to one.

%s9.2 #&#
\subsection{Tightness of error bounds}

In Section~\ref{SECSDP}, we prove that both the $\DD$ and $\MDP_k$
statistics discriminate between $H_0$ and $H_1$ with high probability
as long as $\theta\ge C k\sqrt{\log(p/k)/n}$. The previous
subsection indicates that $\MDP_k$ actually performs better than $\DD
$ and it is pertinent to wonder if detections levels of order smaller
than $\theta\ge C k\sqrt{\log(p/k)/n}$ can be achieved. In this
subsection, we bring numerical evidence that it is not the case and
thus corroborate evidence from Section~\ref{SECCTLB}.

For $\MDP_k$ to be considered a tight (up to constant factor)
approximation of $\lambda_{\max}^k$, it needs to discriminate between
$H_0$ and $H_1$ with high probability as soon as $\theta$ is of the
order $\sqrt{k\log(p/k)/n}$, which is\vspace*{1pt} the minimax optimal detection
level that is also achieved by $\lambda_{\max}^k$. This behavior can
be illustrated by showing a phase transition for the probability of
error in the testing problem, as a function of $\theta$, for different
choices of $(p,n,k)$. More precisely, if $\MDP_k$ were a tight
approximation of $\lambda_{\max}^k$, there should exist a critical
value $\theta_{\mathrm{crit}}$ and a constant $C_{\mathrm{crit}}$,
such that $\theta> \theta_{\mathrm{crit}}=C_{\mathrm{crit}} \sqrt{
k \log(p/k)/n}$, the probability of type II error is close to 0.
Moreover, $C_{\mathrm{crit}}$ should not depend on $(p,n,k)$. Our
numerical results show that this is not the case. Instead, as predicted
by the analysis of Section~\ref{SECSDP}, our experiments point to
$\theta_{\mathrm{crit}}$ of order $k\sqrt{ \log(p/k)/n}$.

In order to substantiate such effects, we use a reciprocal setting. For
fixed $\theta= 1$, fixed probability of type I error (test level) and
several choices of parameters $(p,k)$, we exhibit a phase transition
for the probability of type II error $P_{\mathrm{II}}(\cdot)$ as a
function of the optimal and suboptimal scalings, defined respectively by
\[
\eta^* = \frac{k}{n} \log\biggl(\frac{p}k \biggr) \quad\mbox{and}\quad
\eta^\circ= \frac{k^2}{n} \log\biggl(\frac{p}k \biggr).
\]
If $\eta\in\{\eta^*, \eta^\circ\}$ is the correct scaling, there
should exist a critical value $\eta_{\mathrm{crit}}$, independent of
$(p,n,k)$, such that one of the following two scenarios hold. On the
one hand, if $\MDP_k$ actually exhibits optimal rates, that is, if
$\eta^*$ is the correct scaling, then $\eta^* \mapsto
P_{\mathrm{II}}(\eta^*)$ should have a sharp transition from $0$ to $1$
around
$\eta_{\mathrm{crit}}$ for all choices of parameters $(p,k)$. On the
other hand, if the correct scaling for $\MDP_k$ is~$\eta^\circ$,
then it is the function $\eta^\circ\mapsto P_{\mathrm{II}}(\eta
^\circ)$ that has a sharp transition around $\eta_{\mathrm{crit}}$
for all choices of parameters $(p,k)$.

We simulate $N=1400$ samples of $n$ independent random variables
$X^0_1,\ldots,\allowbreak X^0_n \sim\mathcal{N}(0,I_p)$. It yields $\Sh_1^0,
\ldots, \Sh_N^0$ that are drawn under $H_0$, and used to estimate the
quantiles $q_{0.01},q_{0.05}$ at $1\%$ and $5\%$ for the $\MDP_k$
statistic. The same process is repeated under $H_1$ to estimate the
probability of type II error $\Pro_{H_1}(\MDP_k(\Sh)>q_{\alpha})$.
To that end, we simulate $X^1_1,\ldots,X^1_n \sim\mathcal{N}(0,I_p +
\theta vv^\top)$, for random unit vectors $v$ supported on $S=\{1,\ldots
, k\}$. The restriction of $v$ to $S$ is distributed uniformly
on the unit sphere of dimension $k$. To display a one-dimensional
dependence, $k$ is chosen equal to the integer part of $\sqrt{p}$.

%f2 #&#
%
\begin{figure}

\includegraphics{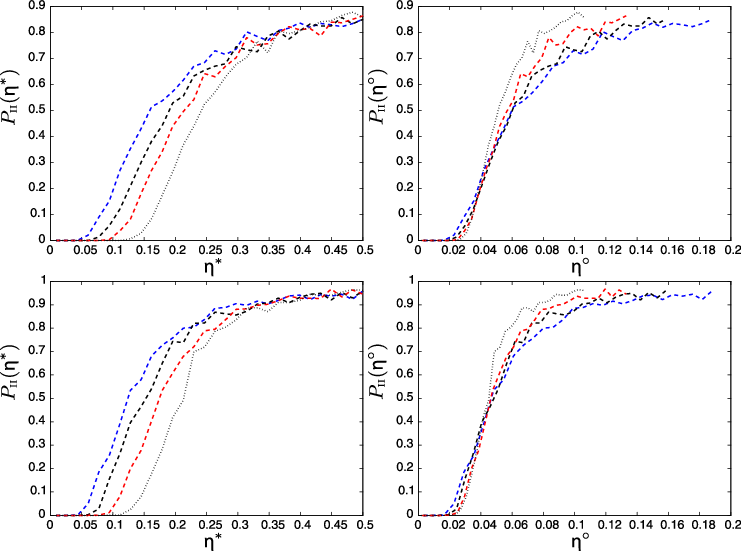}

\caption{Probability of type \textup{II} error $ P_{\mathrm{II}}(\eta)$ with
scalings $\eta=\eta^*$ (left) and $\eta=\eta^\circ$ (right) for
$p=\{50,100,200,500\}$, $k=\lfloor\sqrt p \rfloor$, $N=1400$.
Test levels are $\alpha=5\%$ (top) and $\alpha=1\%$ (bottom).}
\label{FIGscale}\vspace*{-2pt}
\end{figure}

Figure~\ref{FIGscale} compares the behavior of the functions $\eta^*
\mapsto P_{\mathrm{II}}(\eta^*)$ and $\eta^\circ\mapsto
P_{\mathrm{II}}(\eta^\circ)$. It clearly demonstrates the presence
of a
critical level $\eta_{\mathrm{crit}}\simeq0.1$ independent of
$(p,n,k)$. The concomitance of the right curves for different choices
of $(p,n,k)$ indicates that $\eta^\circ$ rather than $\eta^*$ is the
correct scaling factor for the $\MDP_k$ statistic. This confirms the
results of Section~\ref{SECCTLB} and the existence of a statistical
price to pay for computational efficiency.

\begin{appendix}\label{app}
%s10 #&#
\section*{Appendix: Technical lemmas}

We gather in this Appendix various useful concentration inequalities.
The first Lemma is due to Laurent and Massart.

%le10.1 #&#
%
\begin{lem}[{[\citet{LauMas00}, Lemma 1]}]
\label{LEMLM}
Let $Z_1,\ldots,Z_n \sim\cN(0,1)$ be i.i.d. ce random variables, and
define $Y=\frac1n \sum_{i=1}^n Z_i^2 -1$. Then the two following tail
bounds hold for any $t>0$:
\[
\mathbf{P} \biggl(Y \leq-2\sqrt{\frac tn} \biggr) \leq e^{-t},\qquad
\mathbf{P} \biggl(Y \geq2\sqrt{\frac tn} + 2 \frac tn \biggr) \leq
e^{-t}.
\]
\end{lem}

This second lemma generalizes the previous one to sums of squares of
sub-Gaussian random variables.

%le10.2 #&#
%
\begin{lem}
\label{LEMLMsubexp}
Let $G_1,\ldots,G_n$ be i.i.d. standard sub-Gaussian centered random
variables. It holds
\[
\Pro\Biggl( \Biggl|\frac1n \sum_{i=1}^n
\bigl(G_i^2 - \E\bigl[G_i^2
\bigr] \bigr) \Biggr| > 2e \biggl( 64\frac tn + 32 \sqrt{\frac tn} \biggr)
\Biggr) \leq2
e^{-t}.\vadjust{\goodbreak}
\]
\end{lem}
\begin{pf}
Using a Chernoff bound and integrating the tails yields that $\E
[ |G|^p ]^{1/p} \leq2 \sqrt{p}$, for any integer $p\ge0$.
It follows from these bounds, by a series expansion, that
%
%e10.1 #&#
%
\setcounter{equation}{0}
\renewcommand{\theequation}{A.\arabic{equation}}
\begin{equation}
\label{EQMGF2} \E\bigl[e^{t (G^2 - \E[G^2] )} \bigr] \le\exp\bigl(512
e^2 t^2 \bigr)\qquad\mbox{for } 0<t<1/(32e).
\end{equation}
For any $u\in\R^n$, define
\[
S_n = \sum_{i=1}^n
u_i \bigl(G_i^2 - \E\bigl[G_i^2
\bigr] \bigr).
\]
By a Chernoff bound, using equation (\ref{EQMGF2}), it holds for all $t>0$,
\[
\Pro(S_n \ge t) \le\exp\biggl(-\min\biggl(\frac{t^2}{2048 e^2
|u|_2^2},
\frac{t}{64 e |u|_\infty} \biggr) \biggr).
\]
This implies our final result.
\end{pf}

%le10.3 #&#
%
\begin{lem}
\label{LEMabsnorm}
Let $Z_1,\ldots, Z_n$ be i.i.d. $\cN(0,1)$ random variables and
define $
Y=\sum_{i=1}^n |Z_i|$.
Then, for any $t>0$, it holds
\[
\Pro(Y -\E Y< -t) \le e^{-{t^2}/({2n})}.
\]
\end{lem}
\begin{pf}
Using a Chernoff bound, observe first that for any $s>0$,
we have
\begin{eqnarray*}
\Pro(Y -\E Y < -t)&=&\Pro(\E Y-Y> t ) \le e^{-st} \E
\bigl[e^{s(\E
Y-Y)} \bigr]
\\
&=& e^{-st}\prod_{i=1}^n \E
\bigl[e^{-s(\E|Z_i|- |Z_i|)} \bigr].
\end{eqnarray*}
Moreover,
\[
\E\bigl[e^{-s(\E|Z_i|- |Z_i|)} \bigr]\le2e^{-s\E|Z_i|}\E\bigl[e^{sZ_i}
\bigr]=2e^{-s\sqrt{\pi/2}}e^{s^2/2}.
\]
The above two displays yield
\begin{eqnarray*}
\Pro(Y -\E Y < -t)&\le&2^n\inf_{s>0} \exp
\biggl(-st-ns\sqrt{\pi/2}+ n\frac{s^2}{2} \biggr)
\\
&=&2^n\exp\biggl(-\frac{(t+n\sqrt{\pi/2})^2}{2n} \biggr)\le e^{-
{t^2}/({2n})}.
\end{eqnarray*}
\upqed\end{pf}
\end{appendix}


%suskaldyti doi

% imsref loaded by lrinkeviciute, 2013-06-20 13:29:07
% imsref loaded by lrinkeviciute, 2013-06-21 09:48:09
%

\printaddresses

\end{document}